\documentclass{book}
\usepackage[contrib,lang=british]{ems-book} 

\usepackage{bbm}
 \usepackage{tikz, tikz-cd}

 \usetikzlibrary{decorations.markings}
\usetikzlibrary{shapes.geometric, patterns, angles, quotes}
\usetikzlibrary{arrows}
\usetikzlibrary{decorations}
\usetikzlibrary{shapes}
\usetikzlibrary{decorations.pathreplacing,decorations.markings} 

\usetikzlibrary{calc,intersections,through,backgrounds}

\usetikzlibrary{shapes.geometric,positioning}
\usetikzlibrary{arrows,decorations.pathmorphing,decorations.pathreplacing}

\usetikzlibrary{hobby}
\usepackage{adjustbox}

\tikzset{
    clip even odd rule/.code={\pgfseteorule}, 
    invclip/.style={
        clip,insert path=
            [clip even odd rule]{
                [reset cm](-\maxdimen,-\maxdimen)rectangle(\maxdimen,\maxdimen)
            }
    }
}

\newcommand{\cA}{\mathcal{A}}
\newcommand{\cB}{\mathcal{B}}

\newcommand{\cF}{\mathcal{F}}

\newcommand{\sEnd}{\rm\underline{End}} 
\newcommand{\sHom}{{\rm\underline{Hom}}} 

\newcommand{\bA}{\mathbb{A}}
\newcommand{\bB}{\mathbb{B}}

\newcommand{\punct}{\mathscr{P}}

\newcommand{\rmod}[1]{{\rm mod \,}#1} 
\newcommand{\smod}[1]{{\rm \underline{mod} \,}#1} 

\newcommand{\soc}{{\operatorname{soc}}} 
\newcommand{\rad}{{\operatorname{rad}}} 

\DeclareMathOperator{\Hom}{Hom}
\DeclareMathOperator*{\proj}{proj}
\DeclareMathOperator{\End}{End}

\DeclareMathOperator{\marked}{\mathcal{M}}
\DeclareMathOperator{\DA}{\mathsf{D}^b(A)}
\DeclareMathOperator{\DB}{\mathsf{D}^b(B)}

\DeclareMathOperator{\Deck}{\mathsf{Deck}}
\DeclareMathOperator{\H^0}{H^0}
\DeclareMathOperator{\HH}{H}
\DeclareMathOperator{\triv}{T}

\DeclareMathOperator{\Tw}{Tw}
\DeclareMathOperator{\Fuk}{Fuk}
\DeclareMathOperator{\MCG}{MCG}
\DeclareMathOperator{\Arf}{Arf}
\newcommand{\oInt}{\,\vv{\cap}\,}
\newcommand{\Ob}[1]{\operatorname{Ob}(#1)}
\newcommand{\Kb}[1]{\operatorname{K}^b(\proj #1)}

\newcommand{\seq}{\mathbbm{m}}

\newcommand{\BV}[1]{\operatorname{V}(#1)} 
\newcommand{\BE}[1]{\operatorname{E}(#1)} 
\newcommand{\BF}[1]{\operatorname{F}(#1)} 
\newcommand{\Out}[1]{\operatorname{Out^0}(#1)} 
\newcommand{\HHH}[1]{\operatorname{HH^1}(#1)} 
\newtheorem{thm}{Theorem}
\newtheorem{lem}{Lemma}

\begin{document}
\mainmatter

\title{Derived equivalences of Brauer graph algebras}
\titlemark{Derived equivalences of Brauer graph algebras}

\emsauthor{1}{Alexandra Zvonareva}{A.~Zvonareva}


\emsaffil{1}{Postal address \email{Email}}

\classification[16E35]{16G10} 



\begin{abstract}
The aim of this short survey is to trace back the ingredients going into the derived equivalence classification of Brauer graph algebras and into the proof of the fact that these algebras are closed under derived equivalence.
\end{abstract}

\makecontribtitle


\section{Introduction}

The study of derived categories and derived equivalences is one of the major topics in representation theory of finite dimensional algebras. It happens quite rarely that one can describe a class of algebras closed under derived equivalences and then obtain a complete classification of these algebras up to derived equivalence. Gentle algebras provide one example for which this was carried out. Indeed, by \cite{SchroerZimmermann} the class of gentle algebras is closed under derived equivalences. A complete classification of gentle algebras up to derived equivalence was obtained in \cite{AmiotPlamondonSchroll, OpperDerivedEquivalences} (see also \cite{HaidenKatzarkovKontsevich, LekiliPolishchukGentle}). This short survey is dedicated to another class of algebras for which similar results can be obtained: Brauer graph algebras.

First appearing in modular representation theory as blocks with cyclic \cite{Dade} or dihedral defect group \cite{Donovan}, Brauer graph algebras are very well studied. For example, all finite dimensional indecomposable modules over such algebras are classified  \cite{GelfandPonomarev, DonovanFreislich, ButlerRingel, WaldWaschbusch}, the structure of the Auslander-Reiten components is well understood \cite{ErdmannSkowronsk}, and the Yoneda algebras of Brauer graph algebras are known to be finitely generated \cite{AntipovGeneralov, GreenSchrollSnashallTaillefer}.

Throughout this survey all finite dimensional algebras will be algebras over an algebraically closed field $\Bbbk$. As mentioned earlier, we will concentrate on two results. The first being: 

\begin{thm}[\cite{AntipovZvonareva,ZvonarevaAntipov}]
   The class of Brauer graph algebras is closed under derived equivalence.
Namely, if $A$ is  a Brauer graph algebra and $B$ is an
algebra such that $\mathsf{D}^b(\rmod{A}) \simeq \mathsf{D}^b(\rmod{B})$, then $B$ is Morita equivalent to a Brauer graph algebra.
\end{thm}

The second result is the derived equivalence classification of Brauer graph algebras. Brauer graph algebras can be constructed from a ribbon graph $\Gamma$ and  a multiplicity function $\seq$ which assigns a non-zero natural number to each vertex of this graph. The pair $(\Gamma,\seq)$ is called a Brauer graph. Any ribbon graph $\Gamma$ can be embedded into an oriented surface $\Sigma_\Gamma$ in a minimal way. The surface $\Sigma_\Gamma$ is defined  uniquely up to homeomorphism. It turns out that the surface $\Sigma_\Gamma$ plays a crucial role in the study of derived equivalences between Brauer graph algebras. This was first realized by Antipov \cite{Antipov}, who introduced the necessary derived invariants of Brauer graph algebras and studied the case when $\Sigma_\Gamma$ is a sphere with boundary.  In the general case we have the following result:

\begin{thm}[\cite{OpperZvonareva}]
Let $A$ and  $A'$ be two Brauer graph algebras with Brauer graphs $\Gamma$ and $\Gamma'$. Assume that $A$ and $A'$ each have at least two isomorphism classes of simple modules. Then, $A$ and $A'$ are derived equivalent if and only if all of the following conditions are satisfied.
\begin{enumerate}
\setlength\itemsep{1ex}
	\item $\Gamma$ and $\Gamma'$ have the same number of vertices, edges, and faces, in particular, the surfaces of $A$ and $A'$ are homeomorphic;
    \item \label{Enum2} the multisets of perimeters of faces and the multisets of the multiplicities of vertices of $\Gamma$ and $\Gamma'$ coincide;
    \item \label{Enum3} $\Gamma$ and $\Gamma'$ are either both bipartite or both not bipartite.
\end{enumerate}
\end{thm}

The faces of a Brauer graph $\Gamma$ and their perimeters can be defined using the embedding of $\Gamma$ into $\Sigma_\Gamma$. As one sees from this theorem, derived equivalence classes of Brauer graph algebras can be completely described in terms of the corresponding Brauer graphs. The proof of this theorem does not rely on constructing tilting complexes, instead it uses $A_\infty$-categories and connections to Fukaya categories.

\section{Notation}

 Throughout this paper $\Bbbk$ will denote an algebraically closed field. Let $A$ be a finite-dimensional $\Bbbk$-algebra, the category of finite-dimensional right $A$-modules will be denoted by $\rmod{A}$; the bounded derived category of $\rmod{A}$ will be denoted by $\DA$, and the homotopy category of bounded complexes of finitely generated projective $A$-modules by $\Kb{A}$, the stable module category of $A$ will be denoted by $\smod{A}$, $\sHom$ and $\sEnd$ will denote morphism spaces and endomorphism algebras in this category, the syzygy functor on $\smod A$ will be denoted by $\Omega$. The product of arrows in the path algebra of a quiver is read from  right to left, i.e. a path $3\xleftarrow{\beta}2\xleftarrow{\alpha}1$ coincides with the product $\beta\alpha$. For an arrow $\alpha$ we will denote by $s(\alpha)$ its source and by $t(\alpha)$ its target. 

\section{Brauer graph algebras}\label{sec:BGA}

Brauer graph algebras can be defined from the combinatorial data of a ribbon graph $\Gamma$ and a \textbf{multiplicity function} $\seq: \BV{\Gamma}\rightarrow \mathbb{N}$ from the set of vertices of $\Gamma$ to natural numbers. The values of $\seq$ at the vertices of $\Gamma$ will be called \textbf{multiplicities} and the pair $(\Gamma, \seq)$ will be called a \textbf{Brauer graph}.

Roughly speaking a \textbf{ribbon graph} $\Gamma$ is a graph with a cyclic ordering of edges around each vertex. 
This can be formalized by splitting each edge into a pair of half-edges in such a way that each half-edge is adjacent to exactly one vertex of $\Gamma$ and fixing a cyclic ordering of half-edges around vertices of $\Gamma$. For a half edge $h$, the succeeding (resp. preceding) half-edge in the cyclic ordering will be denoted by $h^+$ (resp. $h^-$) and the edge containing $h$ will be denoted by $\overline{h}$. Alternatively, the data of the cyclic ordering of edges can be given by an embedding 
$\iota: \Gamma \hookrightarrow \Sigma$ of $\Gamma$ into the interior of an oriented surface $\Sigma$. Up to homeomorphism, there exists a unique compact, oriented surface $\Sigma_{\Gamma}$ with boundary and a (non-unique) embedding $\iota:\Gamma \hookrightarrow \Sigma_{\Gamma}$ into the interior of $\Sigma_\Gamma$ such that  $\Sigma_{\Gamma} \setminus \Gamma$ is a union of half-open annuli.
 
To any Brauer graph $(\Gamma, \seq)$ one can associate a quiver $Q_\Gamma$ and an ideal of relations $I_\Gamma$ in the path algebra $\Bbbk Q_\Gamma$ as follows.

\begin{enumerate}[wide, labelwidth=!, labelindent=0pt]
    \item The vertices of $Q_\Gamma$ correspond to the edges of $\Gamma$ and for every half-edge $h \in  \Gamma$, there is an arrow $\alpha_h: \overline{h}\rightarrow \overline{h^+}$. The assignment $\alpha_h \mapsto \alpha_{h^-}$ defines a permutation $\pi_{\Gamma}$ of the arrows of $Q_\Gamma$ whose cycles are in bijection with $\BV{\Gamma}$. Hence every arrow $\alpha$ defines a closed path 
\begin{displaymath}
C_{\alpha}=\alpha \pi_\Gamma(\alpha) \cdots \pi_\Gamma^l(\alpha),
\end{displaymath} 

\noindent where $l+1$ is the cardinality of the $\pi_\Gamma$-orbit of $\alpha$. 
Every vertex of $Q_\Gamma$ is the starting point of exactly two cycles of the form $C_{\alpha}$. For an arrow $\alpha$ denote by $\seq(C_{\alpha})$ the multiplicity of the vertex corresponding to the cycle $C_{\alpha}$.

\item The ideal $I_\Gamma$ is generated by the following set of relations:

\begin{enumerate}[wide, labelwidth=!, labelindent=\parindent]
    \item  
    $   \alpha\beta,$  
    \noindent where $\pi_\Gamma(\alpha)\neq \beta$;
     \item  
    $ 
      C_\alpha^{\seq(C_\alpha)} - C_{\beta}^{\seq(C_{\beta})},$ where $\alpha, \beta \in Q_1$ and $t(\alpha)=t(\beta)$, that is $\alpha$ and $\beta$ end at the same edge of $\Gamma$.
\end{enumerate}

\end{enumerate}

 The resulting finite dimensional $\Bbbk$-algebra $\Bbbk Q_{\Gamma} / I_{\Gamma}$ will be denoted by $B_{\Gamma}$. A $\Bbbk$-algebra $B$ is called a \textbf{Brauer graph algebra} if there exists a Brauer graph $(\Gamma, \seq)$ such that $B \cong B_{\Gamma}$ as $\Bbbk$-algebras. For simplicity we will consider only connected Brauer graphs. 
\begin{figure}[H]
\begin{tikzpicture}
\node (v2) at (-3,1.5) {};
\node (v1) at (-3.5,2.5) {};
\node (v3) at (-3.5,0.5) {};
\node (v4) at (-1,1.5) {};
\node at (-0.5,2.5) {};
\node (v5) at (-0.5,2.5) {};
\node (v6) at (-0.5,0.5) {};
\draw  (-3.5,2.5) -- (-3,1.5);
\draw  (-3.5,0.5) edge (-3,1.5);
\draw  (-3,1.5) edge (-1,1.5);
\draw  (-0.5,2.5) edge (-1,1.5);
\draw  (-0.5,0.5) edge (-1,1.5);
\node at (-2,1.5) {};
\node (v8) at (-2,1.5) {};
\node (v9) at (-3.25,2) {};
\node (v7) at (-3.25,1) {};
\node (v10) at (-0.75,1) {};
\node (v11) at (-0.75,2) {};
\draw  (v7) edge[-stealth,bend right] (v8);
\draw  (v8) edge[-stealth,bend right] (v9);
\draw  (v9) edge[-stealth,bend right] (v7);
\draw  (v8) edge[-stealth,bend right] (v10);
\draw  (v10) edge[-stealth,bend right] (v11);
\draw  (v11) edge[-stealth,bend right] (v8);
\draw  (v9) edge[in=160,out=60,loop, looseness=20,-stealth] (v9);
\draw  (v7) edge[in=300,out=200,loop, looseness=20,-stealth] (v7);
\draw  (v10) edge[in=350,out=250,loop, looseness=20,-stealth] (v10);
\draw  (v11) edge[in=110,out=10,loop, looseness=20,-stealth] (v11);
\node at (-2.5,0.85) {$\alpha$};
\node at (-3.6,1.5) {$\beta$};
\node at (-2.5,2.15) {$\gamma$};
\node at (-1.5,2.15) {$\delta$};
\node at (-0.45,1.5) {$\epsilon$};
\node at (-1.5,0.8) {$\zeta$};
\node at (2,2) {$\alpha\beta\gamma=\delta\epsilon\zeta$};
\node at (2,2.5) {$\gamma\delta=0$};
\end{tikzpicture}
    \caption{An example of a Brauer graph $(\Gamma,\seq)$ with $\seq(v)=1$ for all vertices $v\in\Gamma$ together with the quiver $Q_\Gamma$, and one relation of type (a) and one relation of type (b). Note that there are many more relations in $B_\Gamma$.}
    \label{FigureLineFieldRibbonGraph}
\end{figure}
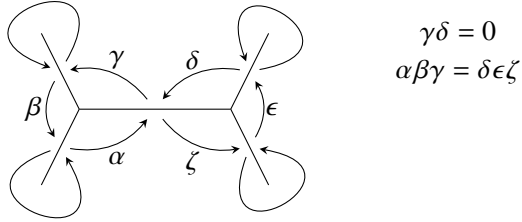
Brauer graph algebras are an example of special biserial algebras, a class of algebras studied quite extensively in representation theory. In fact over an algebraically closed field the class of Brauer graph algebras coincides with the class of symmetric special biserial algebras \cite{SchrollTrivialExtension}. An algebra $A$ is called \textbf{special biserial} if $A$ is Morita equivalent to an algebra $\Bbbk Q / I$ given as a path algebra of a quiver with relations and the following two conditions hold for the pair $(Q,I)$:
\begin{enumerate}[wide, labelwidth=!, labelindent=0pt]
    \item For each vertex $i \in Q$ the number of outgoing arrows at $i$ are less than or equal to $2$; dually, for each vertex $i \in Q$ the number of incoming arrows at $i$ are less than or equal to $2$;
    \item For each arrow $\alpha \in Q$, there is at most one arrow $\beta \in Q$ that satisfies
$\alpha\beta\neq 0$; dually, for each arrow $\alpha \in Q$, there is at most one arrow $\gamma \in Q$ that satisfies
$\gamma\alpha\neq 0$.
\end{enumerate}

Special biserial algebras are tame and their representation theory is very well understood. There is a classification of all indecomposable modules over such algebras in term of strings and bands \cite{WaldWaschbusch}, morphism between indecomposable modules are described in \cite{krause1991maps}, irreducible morphisms and Auslander-Reiten sequences are described in \cite{WaldWaschbusch, ButlerRingel}. The components of the Auslander-Reiten quiver are described in \cite{ErdmannSkowronsk} in the case when the algebra is additionally self-injective (see also \cite{Duffield} for the case of Brauer graph algebras). 

\section{Brauer tree algebras}

Brauer graph algebras of finite representation type can be characterised in terms of their Brauer graphs $(\Gamma,\seq)$ as follows: $\Gamma$ is a tree and $\seq(v)\neq 1$ for at most one vertex $v\in\Gamma$. In that case the vertex $v$ is usually called \textbf{exceptional} and such algebras are called \textbf{Brauer tree algebras}. Note that all non semisimple blocks of group algebras with cyclic defect group are Brauer tree algebras and these are precisely the blocks of finite representation type.

A Brauer graph $(\Gamma,\seq)$ will be called a \textbf{Brauer star} if it has only one vertex $v$ of valency greater than $1$ and the exceptional vertex coincides with $v$. Brauer star algebras $B_{(S,m)}$ for different values $m=\seq(v)$ of the multiplicity of the exceptional vertex constitute the class of basic symmetric Nakayama algebras, i.e. basic symmetric algebras such that every indecomposable module has a unique composition series.

\begin{figure}[h]
\begin{tikzpicture}
\node (v12) at (-3.5,-2) {};
\node (v13) at (-2,-2) {};
\node (v14) at (-3.5,-0.5) {};
\node (v15) at (-5,-2) {};
\node (v16) at (-3.5,-3.5) {};
\draw  (-3.5,-2) edge (v13);
\draw  (-3.5,-2) edge (v14);
\draw  (-3.5,-2) edge (v15);
\draw  (-3.5,-2) edge (v16);
\node (v17) at (-2.5,-2) {};
\node (v18) at (-3.5,-1) {};
\node (v19) at (-4.5,-2) {};
\node (v20) at (-3.5,-3) {};
\draw  (v17) edge[-stealth,bend right] (v18);
\draw  (v18) edge[-stealth,bend right] (v19);
\draw  (v19) edge[-stealth,bend right] (v20);
\draw  (v20) edge[-stealth,bend right] (v17);
\node at (-2.7,-1.2) {$\alpha$};
\node at (-4.3,-1.2) {$\alpha$};
\node at (-4.3,-2.8) {$\alpha$};
\node at (-2.7,-2.8) {$\alpha$};
\node at (-3.3,-1.85) {$m$};
\end{tikzpicture}
    \caption{An example of a Brauer star algebra $B_{(S,m)}$ with $4$ edges. The exceptional vertex is always in the middle. Depending on $m$ all relations are of the form $\alpha^{4m+1}=0$. The loops corresponding to the leaves are omitted, since each loop like that can be expressed as $\alpha^{4m}$. The relations $\alpha^{4m+1}=0$ replace the standard relations.}
    \label{FigureLineFieldRibbonGraph}
\end{figure}
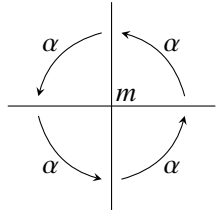

Brauer tree algebras are closed under derived equivalence. Namely, if $A$ is a Brauer tree algebra and $B$ is another finite-dimensional algebra such that $\DA\simeq \DB$, then $B$ is Morita equivalent to a Brauer tree algebra. Let's assume for simplicity that $B_{(S,m)}\not\simeq \Bbbk$ and $B_{(S,m)}\not\simeq \Bbbk [x]/(x^2)$. Note that, since $\Bbbk [x]/(x^2)$ is local and commutative, any algebra derived equivalent to $\Bbbk [x]/(x^2)$ is Morita equivalent to $\Bbbk [x]/(x^2)$ \cite[Corollary 2.13]{RouquierZimmermann}. The fact that Brauer tree algebras are closed under derived equivalence follows from the work of Riedtmann, Gabriel-Riedtmann and Rickard. In \cite{riedtmannalgebren}, Riedtmann in particular shows that any connected algebra $B$ such that $\smod B \simeq \smod B_{(S,m)}$  is Morita equivalent to a Brauer tree algebra. Conversely, in  \cite{gabrielRiedtmann}, Gabriel and Riedtmann show that any Brauer tree algebra is stably equivalent to a symmetric Nakayama algebra. Now, by the work of Rickard, algebras derived equivalent to symmetric algebras are symmetric \cite[Corollary 5.3]{rickardderived}; and derived equivalence for self-injective and, in particular, symmetric algebras implies stable equivalence \cite[Corollary 2.2]{RickardDerivedStable}. Combining these results we get that Brauer tree algebras are closed under derived equivalence. 

The results in \cite{riedtmannalgebren, gabrielRiedtmann} are more general, as the authors consider algebras stably equivalent to all self-injective Nakayama algebras and not only to the symmetric ones. We restricted ourselves to the symmetric case, since this is the case relevant for this note.

Inside the class of Brauer tree algebras the classification up to derived equivalence was established by Rickard \cite[Theorem 4.2]{RickardDerivedStable}, who constructed tilting complexes over Brauer tree algebras whose endomorphism rings are of the form $B_{(S,m)}$.

\begin{thm}[Rickard'89]
    Up to derived equivalence, a Brauer tree algebra is determined by the
number of edges of the Brauer tree and the multiplicity of the exceptional vertex.
\end{thm}

In the subsequent sections we will review the same results for the entire class of Brauer graph algebras. It turns out that Brauer graph algebras are also closed under derived equivalence and their derived equivalence classes can be described. However, the proofs turn out to be more involved.

Derived Picard groups, i.e. the groups of standard autoequivalences of the bounded derived category up to natural isomorphisms, of all self-injective Nakayama algebras were computed in \cite{zvonarevamutations,volkovZvonareva}. In particular, these groups are known for Brauer tree algebras and they give a parametrisation of all tilting complexes over such algebras.

\section{Symmetric stably biserial algebras}

In order to prove that Brauer graph algebras are closed under derived equivalence one can roughly follow the same strategy as in the case of Brauer tree algebras. This strategy was carried out in \cite{Pogorzaly,AntipovZvonareva,ZvonarevaAntipov} in the following steps:

\begin{enumerate}[wide, labelwidth=!, labelindent=0pt]
    \item Describe the class of algebras possibly stably equivalent to Brauer graph algebras, i.e. symmetric stably biserial algebras.
    \item Prove that over a field of characteristic $\neq 2$ symmetric stably biserial algebras are isomorphic to Brauer graph algebras.
    \item Use derived invariants to show that over a field of characteristic $2$ symmetric stably biserial algebras, not isomorphic to Brauer graph algebras, cannot be derived equivalent to such algebras.
\end{enumerate}
 
We will discuss steps (1) and (2) in this section, step (3) will be discussed in more detail in Section \ref{sec:BGAclosed}.
The following result was claimed in \cite[Theorem 7.3]{Pogorzaly}:

\begin{thm}
 Let $A$ be a selfinjective special biserial algebra that is not a Nakayama algebra. If $B$ is stably equivalent to $A$, then $B$ is a  selfinjective special  biserial algebra.
\end{thm}

 As the property of being symmetric is preserved under derived equivalence and derived equivalence of symmetric algebras implies stable equivalence this theorem would imply that the class of symmetric special biserial algebras (the class of algebras coinciding with the class of Brauer graph algebras) is closed under derived equivalence. Unfortunately, in \cite{arikiIijimaPark}, counterexamples for some of the statements of \cite{Pogorzaly} were given.

In \cite[Theorem 1]{AntipovZvonareva} we follow some of the ideas of \cite{Pogorzaly} to reprove the following result:

\begin{thm}\label{thm:sbis}
     Let $A$ be a selfinjective special biserial 
algebra not isomorphic to the Nakayama algebra with $\rad^2=0$. If $B$ is a basic
algebra stably equivalent to $A$, then $B$ is stably biserial.
\end{thm}

 Here a selfinjective algebra $A$ is called \textbf{stably biserial}, if $A\simeq \Bbbk Q/I$ for an admissible ideal of relations $I$, and such that 

 \begin{enumerate}[wide, labelwidth=!, labelindent=0pt]
     \item For each vertex $i \in Q$ the number of outgoing arrows at $i$ are less than or equal to $2$; dually, for each vertex $i \in Q$ the number of incoming arrows at $i$ are less than or equal to $2$;
\item For each arrow $\alpha \in Q$, there is at most one arrow $\beta \in Q$ that satisfies
$\alpha\beta \not\in \soc(A)$; dually for each arrow $\alpha \in Q$, there is at most one arrow $\gamma \in Q$ that satisfies
$\gamma\alpha \not\in \soc(A)$.
 \end{enumerate}
 
The definition above is a simplification of the original definition (see \cite[Proposition 7.5]{arikiIijimaPark}). In the original definition by Pogorza\l{}y, conditions of the form $\alpha\beta \not\in \soc(A)$ are replaced by $\alpha\beta \not\in  \alpha \rad(A) \beta + \soc(A)$. 

The theorem above can be proved using the following idea: if there is an equivalence of categories $\smod B \xrightarrow{F} \smod A$, where
$A$ and $B$ are selfinjective algebras not isomorphic to local Nakayama algebras, then the image of the set of representatives of the
isoclasses of simple $B$-modules is a maximal system of orthogonal stable bricks. An indecomposable $B$-module
$M$ is said to be a \textbf{stable brick} if $\sEnd(M) \simeq \Bbbk$. A family $\{M_i\}_{i\in I}$ of mutually non-isomorphic stable bricks is a \textbf{system of orthogonal stable bricks} if $M_i$ is not of $\tau$-period $1$ for any $i$ and if $\sHom(M_i, M_j) = 0$ for any $i, j \in I$ with $i \neq j$. Finally,
a system of orthogonal stable bricks $\{M_i\}_{i\in I}$ is called \textbf{maximal} if for every
indecomposable $B$-module $N$ that is neither projective nor of $\tau$-period $1$ there
exist $i, j \in I$ such that $\sHom(M_i, N) \neq 0$ and $\sHom(N, M_j) \neq 0$. This notion is closely related to the notion of a simple-minded systems introduced in \cite{KoenigLiu}. Assume $A$ and $B$ to be basic for simplicity and let us fix $A$ to be selfinjective special biserial. Using string combinatorics one can describe all possible maximal systems of orthogonal stable bricks over $A$, studying $\sHom(\tau^{-1} M_j,\oplus_{i \in I}M_i)$ and $\sHom(\tau^{-1} \oplus_{i \in I}M_i, M_j)$ in $\smod A$ and using the Auslander formula one can conclude that in the quiver of $B$ there are at most two incoming and at most two outgoing
arrows at each vertex.

Similarly, employing string combinatorics to study the possible images under $F$ of $B$-modules of the form 
$P/\soc P$ and $\rad P/\soc P$, where $P$ is indecomposable projective, and possible morphisms between these images one concludes that $B$ must be stably biserial. 

On a side note, studying systems of orthogonal stable bricks over stably biserial algebras one can prove one of the Auslander-Reiten conjectures for the class of special biserial algebras:

\begin{thm}[\cite{Pogorzaly,AntipovZvonareva}]
    Let $A, B$ be two finite dimensional
algebras such that $\smod A \simeq \smod B$ and $A$ is special biserial. Then the number of
isomorphism classes of non-projective simple modules over $A$ and $B$ coincides.
\end{thm}

Similarly to the case of symmetric special biserial algebras, symmetric stably biserial algebras can be described using Brauer graphs. For that one determines the $\pi_\Gamma$-cycles as in Section \ref{sec:BGA} and performs a series of base changes. 
Consider the following set-up:

\begin{enumerate}[wide, labelwidth=!, labelindent=0pt]
\item A Brauer graph $(\Gamma, \seq)$;
\item The corresponding quiver $Q_\Gamma$, permutation $\pi_\Gamma$, and cycles $C_\alpha$ for every arrow $\alpha\in Q_\Gamma$ as in Section \ref{sec:BGA};
\item A set $\mathcal{L}=\{\alpha_1,\dots,\alpha_d\}$ of loops, such that $\pi_\Gamma(\alpha_j)\neq \alpha_j$ and a set of elements $\{t_{\alpha_1},\dots,t_{\alpha_d}\},$ with $t_{\alpha_j} \in \Bbbk^*$. 
\end{enumerate}

The loops from the set $\mathcal{L}$ will be called \textbf{deformed loops}. 

\begin{thm}[\cite{AntipovZvonareva,ZvonarevaAntipov}]\label{thm:symstbis}
Any symmetric stably biserial algebra has a presentation $A=\Bbbk Q_\Gamma/I$, where $Q_\Gamma$ is a quiver as in the set-up above and the ideal of relations $I$ is generated by

\begin{enumerate}
\item  $\alpha\beta$ for all arrows $\alpha,\beta \in Q$,
$\beta\neq\pi_\Gamma(\alpha)$, $\alpha\notin\mathcal{L}$;
\item $C_\alpha^{\seq(C_\alpha)}-C_\beta^{\seq(C_\beta)}$ for all $\alpha,\beta\in Q_1$ with $t(\alpha)=t(\beta)$;
\item $\alpha^2-t_{\alpha}C_\alpha^{\seq(C_\alpha)}$ for each $\alpha\in \mathcal{L}$;
\item $C_\alpha^{\seq(C_\alpha)}\beta$ for all $\alpha,\beta\in Q_1$.
\end{enumerate}

\noindent Moreover, any symmetric stably biserial algebra over an algebraically
closed field $\Bbbk$ with characteristic $\neq 2$ is isomorphic to an algebra $\Bbbk Q_\Gamma/I$ as above with   $\mathcal{L}=\emptyset$.
\end{thm}

In particular, any symmetric stably biserial algebra over an algebraically
closed field $\Bbbk$ with characteristic $\neq 2$ is isomorphic to a Brauer graph algebra. So, in that case Brauer graph algebras are closed under derived equivalence. The problem of distinguishing symmetric stably biserial algebras with a non-empty set of deformed loops from Brauer graph algebras over fields of characteristic $2$ up to derived equivalence will be discussed in the next sections. For that, we first have to introduce some known derived invariants of such algebras.

\section{Derived invariants of Brauer graph algebras and  symmetric stably biserial algebras}\label{sec:inv}

Before passing to derived invariants let us observe that the Brauer graph of a symmetric stably biserial algebra is invariant under isomorphism of algebras in almost all cases, so two different Brauer graphs cannot correspond to the same algebra (see \cite[Lemma 3.1]{ZvonarevaAntipov}).

\begin{lem}\label{LemBrauer}
Let $A$ be a symmetric stably biserial algebra with a presentation $\Bbbk Q_\Gamma/I$ as in the previous theorem and the associated Brauer graph $\Gamma$. If $\Gamma$ is not a loop with $1$ as the multiplicity of the unique vertex, or an edge with 2 as the multiplicity of both vertices, then $\Gamma$ does not depend on the choice of the presentation $\Bbbk Q_\Gamma/I$.  
\end{lem}

Since both exceptional cases in the previous lemma correspond to local algebras, from now on we will restrict ourselves to the case when the symmetric stably biserial algebra $A$ has at least two non-isomorphic simple modules. The case when $A$ is local is not very interesting from the point of view of derived equivalences, as they all come from Morita equivalences \cite{RouquierZimmermann}.

Recall that any Brauer graph $\Gamma$ can be embedded into a unique surface $\Sigma_\Gamma$. If we cut the surface $\Sigma_\Gamma$ along the edges of $\Gamma$, we get a collection of polygons, each containing exactly one boundary component. This polygons are called \textbf{faces} of $\Gamma$. The number of edges in such polygon is called its \textbf{perimeter}. Note that an edge of $\Gamma$ can occur twice as an edge of a face, e.g. the surface corresponding to a Brauer star $\Gamma$ with $n$ edges is a sphere with one boundary component. If we cut it along the edges of $\Gamma$ we get a $2n$-gon with one boundary component inside, so it has one face of perimeter $2n$.

\begin{thm}[\cite{Antipov,ZvonarevaAntipov}]\label{TheoremDerivedInvariantsBrauerGraphAlgebra}
Let $A$ be a symmetric stably biserial algebra with a Brauer graph $\Gamma$ and with at least two simple modules. The following are invariants of $A$ under derived equivalence of symmetric stably biserial algebras:
\begin{enumerate}
    \item \label{Surfaces} the number of vertices, the number of edges, and the number of faces of $\Gamma$;
    \item \label{NumberFaces} the multiset of perimeters of faces of $\Gamma$;
    \item the multiset of multiplicities of vertices of $\Gamma$;
    
    \item  whether $\Gamma$ is bipartite.
\end{enumerate}
\end{thm}

One can easily see from the theorem that the surfaces $\Sigma_\Gamma$ and $\Sigma_{\Gamma'}$ corresponding to two derived equivalent stably biserial algebras with Brauer graphs $\Gamma$ and $\Gamma'$ are homeomorphic. 

The corresponding results for Brauer graph algebras were proved by Antipov with some minor inaccuracies in \cite{antipov2005grothendieck,antipov2006invariants,Antipov,antipov2009structure}, the proofs for symmetric stably biserial algebras are identical or rely on the
corresponding results for Brauer graph algebras and can be found in \cite{ZvonarevaAntipov}. 

Let us have a closer look at the algebraic interpretation of these combinatorial invariants. Assume that the symmetric stably biserial algebra $A$ has a presentation as in Theorem~\ref{thm:symstbis} with the corresponding Brauer graph $\Gamma$. Note that the loops $\alpha \in Q_\Gamma$ such that $\pi_\Gamma(\alpha)\neq \alpha$ are exactly the loops of $Q_\Gamma$ coming from the faces of $\Gamma$ of perimeter $1$.

\begin{enumerate}[wide, labelwidth=!, labelindent=0pt]
    \item The number of edges of the Brauer graph corresponds to the number of isomorphism classes of simple $A$-modules, so it coincides with the rank of the Grothendieck group of $\DA$.
    \item The description of the center of $A$ allows to conclude that the multiset of multiplicities of vertices different from $1$ and the number of faces of $\Gamma$ of perimeter $1$ is a derived invariant.
    \item In the Auslander-Reiten quiver of $\smod A$
each periodic component is a tube. Moreover, all tubes are either tubes of rank $1$, consisting of band modules, or tubes consisting of string modules, the latter tubes are called \textbf{exceptional}. Exceptional tubes correspond to faces  of $\Gamma$: if a face has an even perimeter $p$, then it produces two tubes of rank $p/2$, which are permuted by $\Omega$ (if $p/2$ is at least $2$); if a face has an odd perimeter, then it produces one tube of rank $p$, which is stable under the action of $\Omega$. This allows to deduce that the number and the perimeters of faces of perimeter at least $3$ is invariant under derived equivalence. The invariance of the number of faces or perimeter $2$ can be deduced from the invariance of the number of arrows in $Q_\Gamma$ (the number of arrows is exactly $2|\BE \Gamma|$) and the invariance of the number of the other faces (note that the quiver $Q_\Gamma$ is not necessarily admissible for $A$). 
\item The rank of the Grothendieck group of $\smod A$ is $|\BE \Gamma|-|\BV \Gamma|+1$ if $\Gamma$ is bipartite and $|\BE \Gamma|-|\BV \Gamma|$, if not. Here $|\BV \Gamma|$, $|\BE \Gamma|$, and $|\BF \Gamma|$ denotes the number of vertices, edges, and faces of $\Gamma$, respectively.
We already know that the number of faces and the number of edges of $\Gamma$ are derived invariants. Since $|\BV \Gamma|-|\BE \Gamma|+|\BF \Gamma|$ is even as the Euler characteristic of the surface $\Sigma_\Gamma$ with the boundary components glued in, we see that the parity of $|\BV \Gamma|$ is a derived invariant, so $|\BV \Gamma|$ and the bipartivity of $\Gamma$ must be derived invariants as well.
\end{enumerate}

\section{Brauer graph algebras are closed under derived equivalence}\label{sec:BGAclosed}

In this section we will discuss an additional derived invariant of symmetric stably biserial algebras which can be used in characteristic $2$ to prove that the number of deformed loops is a derived invariant. This invariant is the identity component $\Out A$ of the group of outer automorphisms of $A$. By the results of Huisgen-Zimmermann and Saor\'in, and Rouquier, for a
finite dimensional algebra $A$ over an algebraically closed field the group $\Out A$ is invariant under derived equivalence as an algebraic group \cite{huisgen2001geometry,rouquier2011automorphismes}. This invariant is not used very often. The only application we know of is the proof of
the fact that the number of arrows in the quiver of a gentle algebra is a derived invariant
\cite{AvellaAlaminosGeiss}. 

In the following theorem we need to exclude one further family of algebras in addition to algebras with a unique isomorphism class of simple modules. Namely, the algebras called \textbf{caterpillar} algebras in \cite{ZvonarevaAntipov}. These algebras have Brauer graphs of one of the following two forms (the cyclic ordering of edges comes from the counter-clock-wise orientation on the plane):

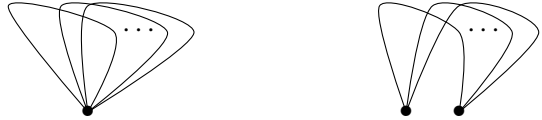
\begin{center}
\begin{tikzpicture}[scale=0.7]

\node (v1) at (0.5,-0.5) {$\bullet$};
\draw  plot[smooth, tension=.7] coordinates {(v1) (2.5,1) (0.5,1.5) (0.5,-0.5)};
\draw  plot[smooth, tension=.7] coordinates {(0.5,-0.5) (2,1) (0,1.5) (0.5,-0.5)};
\node at (1.5,1) {$\cdots$};
\draw  plot[smooth, tension=.7] coordinates {(0.5,-0.5) (1,1) (-1,1.5) (0.5,-0.5)};
\node (v3) at (6.5,-0.5) {$\bullet$};
\node (v2) at (7.5,-0.5) {$\bullet$};
\draw  plot[smooth, tension=.7] coordinates {(v2) (9,1) (7.5,1.5) (v3)};
\draw  plot[smooth, tension=.7] coordinates {(7.5,-0.5) (8.5,1) (7,1.5) (6.5,-0.5)};
\node at (8,1) {$\cdots$};
\draw  plot[smooth, tension=.7] coordinates {(7.5,-0.5) (7.5,1) (6,1.5) (6.5,-0.5)};
\end{tikzpicture}
\end{center}

The quiver in both cases is of the following form which explains the name ($n> 1$):

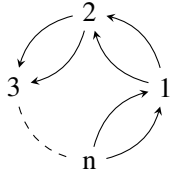
\begin{center}
\begin{tikzpicture}
\node (v17) at (-2.5,-2) {1};
\node (v18) at (-3.5,-1) {2};
\node (v19) at (-4.5,-2) {3};
\node (v20) at (-3.5,-3) {n};
\draw  (v17) edge[-stealth,bend right] (v18);
\draw  (v18) edge[-stealth,bend right] (v19);
\draw  (v19) edge[dashed,bend right] (v20);
\draw  (v20) edge[-stealth,bend right] (v17);

\draw  (v17) edge[-stealth,bend left] (v18);
\draw  (v18) edge[-stealth,bend left] (v19);
\draw  (v20) edge[-stealth,bend left] (v17);

\end{tikzpicture}
\end{center}

Brauer graphs of caterpillar algebras have no faces of perimeter $1$. So any symmetric stably biserial algebra derived equivalent to a caterpillar algebra has no such faces in its Brauer graph as well, and thus must be a Brauer graph algebra, as there are no loops, which can be deformed.  We get that caterpillar algebras cannot be derived equivalent to a symmetric stably biserial algebra, which is not a Brauer graph algebra. 

\begin{thm}[\cite{ZvonarevaAntipov}]\label{TheoremRank} 
Let $\Bbbk$ be an algebraically closed field.  Let $A$ be a symmetric stably biserial algebra over $\Bbbk$  with at least two non-isomorphic simple modules, which is not a caterpillar. Let $\Gamma$ be the Brauer graph of $A$ and let $d$ be the number of deformed loops in $A$ ($d=0$ if  characteristics of $\Bbbk$ is different from $2$). The rank of the maximal
torus of $\Out A$ is $|\BE \Gamma|-|\BV \Gamma|-d+2$.
\end{thm}

Note that maximal tori in $\Out A$ are related to the maximal tori in the Hochschild cohomology Lie algebra $\HHH A$ (see \cite{briggs2023maximal}).

Since $|\BE \Gamma|$ and $|\BV \Gamma|$ are derived invariants of symmetric stably biserial algebras, we immediately get that the number of deformed loops is a derived invariant as well. Indeed, if $A$ and $A'$ are derived equivalent $\Out A\simeq \Out{A'}$ by \cite{huisgen2001geometry,rouquier2011automorphismes}, so the rank of the maximal tori of these two algebraic groups must coincide. Using the fact that over a field of characteristic $2$ Brauer graph algebras are symmetric stably biserial algebras with $0$ deformed loops, we get that they are closed under derived equivalence.

\section{Kauer moves and the case of a sphere}

One of the simplest derived equivalences between Brauer graph algebras is given by tilting mutation. This was first introduced by Kauer long before silting and tilting mutation were studied so intensely \cite{Kauer}. Consider a Brauer graph algebra $A$ as a tilting complex over itself, $A$ can be decomposed into a direct sum of indecomposable projective modules $P_i$, corresponding to the edges of $\Gamma$: $A =\bigoplus\limits_{i\in \BE \Gamma}P_i$. An irreducible silting mutation of $A$ with respect to a fixed indecomposable projective summand $P_j$ is a tilting complex $T=\bigoplus\limits_{i\in \BE \Gamma, i\neq j}P_i \oplus T_j$, since over a symmetric algebra any silting complex is tilting. The complex $T_j$ is computed as the cone of the minimal left approximation of $P_j$ with respect to the additive closure of $\bigoplus\limits_{i\in \BE \Gamma, i\neq j}P_i$. The complex $T_j$ can be easily described depending on the edge $j$. There are three possibilities:
\begin{enumerate}[wide, labelwidth=!, labelindent=0pt]
    \item $j$ is a leaf;
    \item $j$ is a loop encircling a face of perimeter $1$;
    \item $j$ is neither a leaf nor a loop encircling a face of perimeter $1$.
\end{enumerate}    
Let $j_1,j_2$ be the half-edges constituting $j$. In the first two cases denote by $j^+$ the half-edge following $j_1$ or $j_2$ in the cyclic ordering and such that $j^+\neq j_1,j_2$ (there is a unique such half-edge). Depending on the case $T_j$ has the following form: 
\begin{enumerate}[wide, labelwidth=!, labelindent=0pt]
    \item $T_j=P_j\xrightarrow{\alpha}P_{j^+}$, where $\alpha$ is the arrow $j\rightarrow \bar{j^+}$;
    \item $T_j=P_j\xrightarrow{\binom{\alpha}{\alpha\gamma}}P_{j^+} \oplus P_{j^+}$, where $\alpha$ is the arrow $j\rightarrow \bar{j^+}$ and $\gamma$ is the loop $j\rightarrow j$;
    \item $T_j=P_j\xrightarrow{\binom{\alpha}{\beta}}P_{j_1^+} \oplus P_{j_2^+}$, where $\alpha$ is the arrow $j\rightarrow \bar{j_1^+}$ and $\beta$ is the arrow $j\rightarrow \bar{j_2^+}$.
\end{enumerate}

What makes the tilting complexes $T$ particularly nice is that the Brauer graph corresponding to the endomorphism algebra of $T$ can be described as a local geometric operation called the Kauer move. Depending on the case the Kauer move takes the following form:

\begin{figure}[H]
\begin{tikzpicture}[scale=0.83]
    \node (v22) at (-3.5,-6) {};
\node (v15) at (-3.5,-6) {};
\node at (-5.5,-6) {};
\node (v13) at (-5.5,-6) {};
\node (v16) at (-4,-7) {};
\node (v21) at (-3,-7) {};
\node (v14) at (-5,-7) {};
\node (v12) at (-6,-7) {};
\draw  (v12) edge (-5.5,-6);
\draw  (-5.5,-6) edge (v14);
\draw  (-5.5,-6) edge (-3.5,-6) ;
\draw  (-3.5,-6)  edge (v16);
\draw  (-3.5,-6)    edge (v21);

\node (v28) at (-1.5,-6) {};
\node (v29) at (0,-6) {};
\draw  (v28) edge[-stealth]  (v29);
\node (v30) at (1.5,-7) {};
\node (v31) at (2,-6) {};
\node (v32) at (2.5,-7) {};
\node (v33) at (4,-6) {};
\node at (3.5,-7) {};
\node (v34) at (3.5,-7) {};
\node (v35) at (4.5,-7) {};
\draw  (v30) edge (2,-6);
\draw  (2,-6) edge (v32);
\draw  (2,-6)  edge (4,-6);
\draw  (4,-6) edge (v34);
\draw  (4,-6) edge (v35);

\node at (-4.1,-5.5) {\small{$\alpha$}};

\node at (-3.5,-4.65) {\small$j$};

\node (v24) at (-3.5,-5) {};
\draw (-3.5,-6)
 edge (v24);
\node (v23) at (-3.5,-5.5) {};
\node (v25) at (-4,-6) {};
\draw  (v23) edge[-stealth,bend right] (v25);
\node (v26) at (2,-5) {};
\draw  (v26) edge (2,-6)
;
\end{tikzpicture}
    \caption{Local change under the Kauer move from the Brauer graph of $A$ (on the left) to the Brauer graph of $\End_{\Kb{A}} T$ (on the right). Case (1).}
    \label{fig:enter-label}
\end{figure}

\begin{figure}[H]
    
    \begin{tikzpicture}[scale=0.83]
\node (v22) at (-3.5,-6) {};
\node (v15) at (-3.5,-6) {};
\node at (-5.5,-6) {};
\node (v13) at (-5.5,-6) {};
\node (v16) at (-4,-7) {};
\node (v21) at (-3,-7) {};
\node (v14) at (-5,-7) {};
\node (v12) at (-6,-7) {};
\draw  (v12) edge (-5.5,-6);
\draw  (-5.5,-6) edge (v14);
\draw  (-5.5,-6) edge (-3.5,-6) ;
\draw  (-3.5,-6)  edge (v16);
\draw  (-3.5,-6)   edge (v21);

\node (v28) at (-1.5,-6) {};
\node (v29) at (0,-6) {};
\draw  (v28) edge[-stealth]  (v29);
\node (v30) at (1.5,-7) {};
\node (v31) at (2,-6) {};
\node (v32) at (2.5,-7) {};
\node (v33) at (4,-6) {};
\node at (3.5,-7) {};
\node (v34) at (3.5,-7) {};
\node (v35) at (4.5,-7) {};
\draw  (v30) edge (2,-6);
\draw  (2,-6) edge (v32);
\draw  (2,-6)  edge (4,-6);
\draw  (4,-6) edge (v34);
\draw  (4,-6) edge (v35);

\node at (-4.3,-5.55) {\small{$\alpha$}};
\node at (-3.5,-5.15) {\small{$\gamma$}};
\node at (-3.5,-4.3) {\small$j$};
\draw  (v22) edge[loop, looseness=20] (v22);
\node (v23) at (-3.05,-5.5) {};
\node (v24) at (-3.9,-5.5) {};
\node (v25) at (-4.25,-6) {};
\draw  (v23) edge[-stealth,bend right] (v24);
\draw  (v24) edge[-stealth,bend right] (v25);
\draw  (v31) edge[loop, looseness=20] (v31);
\end{tikzpicture}
    \caption{Local change under the Kauer move from the Brauer graph of $A$ (on the left) to the Brauer graph of $\End_{\Kb{A}} T$ (on the right). Case (2).}
    \label{fig:enter-label}
\end{figure}

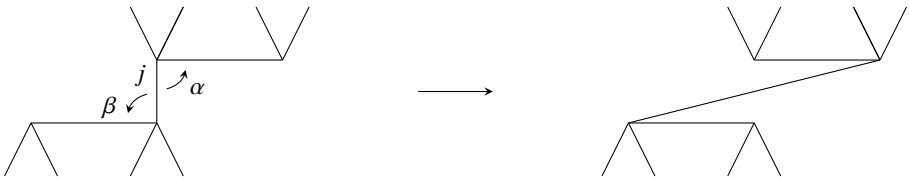
\begin{figure}[H]
    \centering
    \begin{tikzpicture}[scale=0.83]
    \node (v22) at (-3.5,-5) {};
\node (v23) at (-1.5,-5) {};
\node (v24) at (-2,-4) {};
\node (v25) at (-1,-4) {};
\node at (-3.5,-6) {};
\node (v15) at (-3.5,-6) {};
\node at (-5.5,-6) {};
\node (v13) at (-5.5,-6) {};
\node (v16) at (-4,-7) {};
\node (v21) at (-3,-7) {};
\node (v14) at (-5,-7) {};
\node (v12) at (-6,-7) {};
\draw  (v12) edge (-5.5,-6);
\draw  (-5.5,-6) edge (v14);
\draw  (-5.5,-6) edge (-3.5,-6) ;
\draw  (-3.5,-6)  edge (v16);
\draw  (-3.5,-6)  edge (v21);
\draw  (-3.5,-6)  edge (-3.5,-5);
\draw  (-3.5,-5) edge (-1.5,-5);
\draw  (-1.5,-5) edge (v24);
\draw  (-1.5,-5) edge (v25);
\node (v26) at (-3,-4) {};
\node (v27) at (-4,-4) {};
\draw  (-3.5,-5) edge (v26);
\draw  (v26) edge (-3.5,-5);
\draw  (-3.5,-5) edge (v27);
\node (v28) at (0.5,-5.5) {};
\node (v29) at (2,-5.5) {};
\draw  (v28) edge[-stealth]  (v29);
\node (v30) at (3.5,-7) {};
\node (v31) at (4,-6) {};
\node (v32) at (4.5,-7) {};
\node (v33) at (6,-6) {};
\node at (6,-5) {};
\node (v36) at (6,-5) {};
\node at (8,-5) {};
\node (v39) at (8,-5) {};
\node at (5.5,-4) {};
\node (v37) at (5.5,-4) {};
\node (v38) at (6.5,-4) {};
\node at (7.5,-4) {};
\node (v40) at (7.5,-4) {};
\node at (8.5,-4) {};
\node (v41) at (8.5,-4) {};
\node at (5.5,-7) {};
\node (v34) at (5.5,-7) {};
\node (v35) at (6.5,-7) {};
\draw  (v30) edge (4,-6);
\draw  (4,-6) edge (v32);
\draw  (4,-6) edge (6,-6);
\draw  (6,-6) edge (v34);
\draw  (6,-6) edge (v35);
\draw  (6,-5) edge (v37);
\draw  (6,-5) edge (v38);
\draw  (6,-5) edge (8,-5);
\draw  (8,-5) edge (v40);
\draw  (v40) edge (8,-5);
\draw  (v41) edge (8,-5);
\draw  (8,-5) edge (4,-6);
\node (v44) at (-4,-6) {};
\node at (-3.5,-5.5) {};
\node (v42) at (-3.5,-5.5) {};
\node (v43) at (-3,-5) {};
\draw  (v42) edge[-stealth,bend right] (v43);
\draw  (v42) edge[-stealth,bend right] (v44);
\node at (-2.85,-5.45) {\small{$\alpha$}};
\node at (-4.25,-5.75) {\small{$\beta$}};
\node at (-3.75,-5.25) {\small{$j$}};
\end{tikzpicture}
    \caption{Local change under the Kauer move from the Brauer graph of $A$ (on the left) to the Brauer graph of $\End_{\Kb{A}} T$ (on the right). Case (3).}
    \label{fig:enter-label}
\end{figure}

Using Kauer moves Antipov classified up to derived equivalence all Brauer graph  algebras whose Brauer graph can be embedded into a sphere with boundary \cite[Theorem 5.22]{Antipov}. It turns out that it is enough to use the combinatorial invariants discussed in Section \ref{sec:inv}.

\begin{thm}[Antipov'07]
Let $A_1$ and $A_2$ be Brauer graph algebras with Brauer graphs $\Gamma_1$ and $\Gamma_2$ such that $\Sigma_{\Gamma_1}$ and $\Sigma_{\Gamma_2}$ are both homeomorphic to spheres with boundary. Then $A_1$ and $A_2$ are derived equivalent
if and only $\Gamma_1$ and $\Gamma_2$ have the same multisets of perimeters of faces and the same multisets of
multiplicities of vertices.    
\end{thm}

The bipartivity does not appear as one of the invariants in the previous theorem, since  a graph embedded
into a sphere with boundary is bipartite if and only if the perimeters of all its faces are even.

\section{Derived equivalence classification of Brauer graph algebras}

Similar to the case of the sphere in order to classify all Brauer graph algebras up to derived equivalence it is enough to use only the invariants from Section \ref{sec:inv}.

 \begin{thm}[\cite{OpperZvonareva}]\label{IntroTheoremCriterionDerivedEquivalence}
Let $A$ and  $A'$ be two Brauer graph algebras over an algebraically closed field with Brauer graphs $\Gamma$ and $\Gamma'$, assume that $A$ and $A'$ each have at least two isomorphism classes of simple modules. Then, $A$ and $A'$ are derived equivalent if and only if all of the following conditions are satisfied.
\begin{enumerate}
\setlength\itemsep{1ex}
	\item $\Gamma$ and $\Gamma'$ have the same number of vertices, edges and faces, in particular, the surfaces of $A$ and $A'$ are homeomorphic;
    \item \label{Enum2} the multisets of perimeters of faces and the multisets of the multiplicities of vertices of $\Gamma$ and $\Gamma'$ coincide;
    \item \label{Enum3} $\Gamma$ and $\Gamma'$ are either both bipartite or both not bipartite.
\end{enumerate}
\end{thm}

The proof of this theorem is rather involved and uses a lot of ingredients:  $A_\infty$-structures, the connection between Brauer graph algebras and gentle algebras, the connection between gentle algebras and Fukaya categories. Here comes an informal sketch of some of the ideas going into the proof (for details see \cite{OpperZvonareva}).

\subsection{Brauer graph algebras and gentle algebras}

The connection between Brauer graph algebras and gentle algebras was studied by Schroll \cite{SchrollTrivialExtension}. It turns out that some Brauer graph algebras can be obtained as trivial extensions of gentle algebras. Recall that for a finite-dimensional $\Bbbk$-algebra $A$, its \textbf{trivial extension} $\triv(A)$ is the algebra on the vector space $A\oplus \Hom_\Bbbk(A,\Bbbk)$ whose multiplication is given by the formula $(a,f)(b,g)\coloneqq (ab, ag+fb)$, where $ag$ and $fb$ are defined using the  action of $A$ on the $A\text{-}A$-bimodule $\Hom_{\Bbbk}(A,\Bbbk)$. 

In order to pass from a Brauer graph algebra $B$ with a Brauer graph $\Gamma$ (such that $\seq(v)=1$ for every vertex $v\in \Gamma$) to a gentle algebra $A$ one can perform the following operation: fix an arrow $\alpha_v$ in every $\pi_\Gamma$-cycle of $Q_\Gamma$ (note that these cycles are parameterised by the vertices $v\in \Gamma$), consider the ideal $I$, generated by the arrows $\{\alpha_v\}_{v\in \Gamma}$, consider the quotient $A=B/I$. The algebra $A$ is gentle and every basic gentle algebra can be obtained in this way. So we will use this construction instead of the definition, which can be given in a similar way to the definition of a special biserial algebra.  

The connection between gentle algebras and Brauer graph algebras can be summarised as follows:

\begin{thm}[Schroll]
    Let $A$ be a gentle algebra over $\Bbbk$. Then $\triv(A)\simeq B$, where $B$ is a Brauer graph algebras with multiplicity of all vertices $1$. Conversely, let $B$ be a Brauer graph algebras with multiplicity of all vertices $1$, then $B\simeq \triv(A)$ for a gentle algebra $A$. 
\end{thm}

Starting from a Brauer graph algebra (with multiplicities of all vertices $1$), taking the quotient  $A=B/I$, and then taking the trivial extension reconstructs the initial Brauer graph algebra $B\simeq \triv(B/I)$. Gentle algebras can be constructed from graphs on surfaces in a similar way to Brauer graph algebras, the difference is that the vertices of the graph are now supposed to lie on the boundary of the surface.  Pictorially the process of passing from a Brauer graph algebra to a gentle algebra and back can be expressed as follows:

\begin{figure}[H]
    \centering
    
\begin{tikzpicture}
[scale=0.86]

\node (v7) at (0.5,-2.6) {};
\node (v1) at (0.5,-3.05) {};
\draw  (v7)  ellipse (2.5 and 2.5);
\draw  (v1) ellipse (0.5 and 0.5);

\fill[gray, opacity=0.1, even odd rule]  (v7)  ellipse (2.5 and 2.5) (v1) ellipse (0.5 and 0.5);

\node (v2) at (0.5,-1) {$\bullet$};
\node (v4) at (1.5,-4) {$\bullet$};
\node (v3) at (-0.5,-4) {$\bullet$};
\draw  (0.5,-1) edge[bend right] (-0.5,-4);
\draw  (1.5,-4) edge[bend right] (0.5,-1);
\draw  (-0.5,-4) edge (1.5,-4);
\node (v5) at (0.5,-4) {};
\node (v6) at (-0.5,-2.5) {};
\node (v7) at (1.5,-2.5) {};
\draw  (v5) edge[-stealth,bend left] (v6);
\draw  (v5) edge[stealth-,out=240,in=180, looseness=2] (v6); 
\draw  (v5) edge[stealth-,bend right] (v7);
\draw  (v7) edge[stealth-,out=0,in=300, looseness=2] (v5);
\draw  (v7) edge[stealth-,bend right] (v6);
\draw  (v6) edge[stealth-,out=120,in=60, looseness=3] (v7);

\node (v17) at (7.6,-2.6) {};
\node (v11) at (7.6,-3.05) {};
\draw  (v17)  ellipse (2.5 and 2.5);
\draw  (v11) ellipse (0.5 and 0.5);
\fill[gray, opacity=0.1, even odd rule]  (v17) ellipse (2.5 and 2.5) (v11) ellipse (0.5 and 0.5);

\node (v2) at (7.6,-0.1) {$\bullet$};
\node (v4) at (8.85,-4.75) {$\bullet$};
\node (v3) at (7.15,-3.3) {$\bullet$};
\draw  (7.6,-0.1) edge[bend right] (7.15,-3.3);
\draw  (8.85,-4.75) edge[bend right] (7.6,-0.1);
\draw  (7.15,-3.3) edge[bend right] (8.85,-4.75);
\node (v5) at (7.5,-4) {};
\node (v6) at (6.85,-2.2) {};
\node (v7) at (8.9,-2.3) {};

\draw  (v5) edge[stealth-,out=180,in=220, looseness=1] (v6); 
\draw  (v5) edge[stealth-,bend right] (v7);

\draw  (v7) edge[stealth-,bend right] (v6);

\node (v8) at (3.5,-3) {};
\node (v9) at (4.5,-3) {};
\draw  (v8) edge[stealth-stealth] (v9);
\node at (0.5,-0.5) {\small{$\alpha$}};
\node at (2.1,-4.1) {\small{$\beta$}};
\node at (-0.3,-3.75) {\small{$\gamma$}};

\end{tikzpicture}

    \caption{The Brauer graph $\Gamma$ embedded into $\Sigma_\Gamma$ and the quiver $Q_\Gamma$ on the left. The ideal $I$ is generated by the arrows $\alpha,\beta,\gamma$. The quiver of the corresponding gentle algebra $B_\Gamma/I$ on the right. The vertices of the graph $\Gamma$ are dragged to the boundary components of $\Sigma_\Gamma$, the component for the vertex $v$ is chosen based of the arrow $\alpha_v$.}
    \label{gentle drag}
\end{figure}
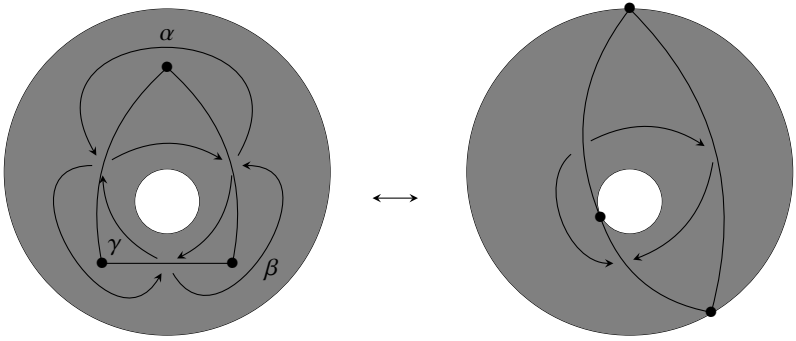

\subsection{$A_\infty$-categories and trivial extensions}

$A_\infty$-algebras were introduced by Stasheff in the sixties in order to study topological spaces. For convenience we will work with $A_\infty$-categories. The sign conventions we use follow some of the literature on Fukaya categories (the disadvantage of this sign convention is that the composition $\mu^2$ is not associative even if $\mu^1$ is zero; one can fix this by twisting $\mu^2$ by an appropriate sign).

An \textbf{$A_\infty$-category}  $\bA$ consists of a set of objects $\Ob{\mathbb{A}}$,
a $\mathbb{Z}$-graded $\Bbbk$-vector space $\Hom_\mathbb{A}(X,Y)$  for each pair of objects $X,Y \in \Ob{\mathbb{A}}$, and graded $\Bbbk$-linear maps
$$
\begin{tikzcd} \mu^n:=\mu^n_\bA : \Hom_\mathbb{A}(X_{n-1},X_n) \otimes \cdots \otimes \Hom_\mathbb{A}(X_0,X_1) \arrow{r} & \Hom_\mathbb{A}(X_0,X_n) \end{tikzcd}$$
 of degree $2 - n$ for each $n \geq 1$, satisfying the \textbf{$A_\infty$-relations}:
\[
   \sum_{k,m}(-1)^{||a_k||+\cdots +||a_1||}\mu^{n-m+1}(a_n,..., a_{k+m+1},\mu^{m}(a_{k+m},..., a_{k+1}), a_k,..., a_1)=0,
\]
 for each sequence of composable elements $a_n,\dots,a_2,a_1$, where $||a_i||=|a_i|-1$.

We will consider $A_{\infty}$-categories which are \textbf{strictly unital}, i.e. such that for every object $X \in \Ob{\mathbb{A}}$ there is an element $1_X \in \Hom_{\mathbb{A}}^0(X, X)$ such that 
$\mu^1(1_X)=0,$ and  $\mu^2(a, 1_X) = (-1)^{|a|}\mu^2(1_Y , a) = a$, for all $a \in \Hom_{\mathbb{A}}(X, Y)$, and $\mu^n(\dots, 1_X, \dots) = 0$ for all $n \geq 3$.

To an $A_{\infty}$-category one can associate another $A_{\infty}$-category $\Tw\bA$ -  the \textbf{category of twisted complexes} over $\bA$. The homotopy category $\H^0(\Tw\bA)$ of $\Tw\bA$ is triangulated and can be thought of as a generalization of the homotopy category of projectives over an ordinary algebra. Indeed, to a basic finite-dimensional algebra $A$ one can associate  an  $A_\infty$-category $\bA$ with objects corresponding to the indecomposable projective $A$-modules, morphisms corresponding to morphisms between them, and $\mu_{\mathbb{A}}^2$ corresponding to composition. The category $\H^0(\Tw\bA)$ is then equivalent to the category $\Kb{A}$. In particular, if we want to produce equivalences between derived categories of two ordinary algebras $\DA\simeq\DB$, or equivalently between the homotopy categories of projectives over these algebras $\Kb{A} \simeq \Kb{B}$, we can work with the categories $\H^0(\Tw\bA)$ and $\H^0(\Tw\bB)$ instead.

Trivial extensions of $A_\infty$-categories were considered in the literature \cite{SeidelSub}. For an $A_\infty$-category $\mathbb{A}$ its \textbf{trivial extension} is an $A_{\infty}$-category 
$\triv(\mathbb{A})$ with the same objects as $\mathbb{A}$. For $X,Y \in \triv(\mathbb{A})$ we have
$\Hom_{\triv(\mathbb{A})}(X,Y) \coloneqq \Hom_{\mathbb{A}}(X,Y) \oplus \big(\mathbb{D} \Hom_{\mathbb{A}}(Y,X)\big),$ where $\mathbb{D}$ denotes the dual graded vector space.

The multi-linear maps $\mu_{\triv(\mathbb{A})}^r$ are uniquely determined as follows:
\begin{enumerate}[wide, labelwidth=!, labelindent=0pt]
\item For $a_r, \dots, a_1$ with $a_j \in \Hom_{\mathbb{A}}(X_{j-1}, X_{j})$,
$$
 \mu_{\triv(\mathbb{A})}(a_r, \dots, a_1) = \mu_{\mathbb{A}}(a_r, \dots, a_1);
$$
\item For $a_{r},\dots, a_{i+1}, a_{i-1},\dots, a_1$ with $a_j \in \Hom_{\mathbb{A}}(X_{j-1}, X_{j})$ and $f \in \mathbb{D}\Hom_{\mathbb{A}}(X_i, X_{i-1})$,
$$
 \mu_{\triv(\mathbb{A})}(a_r, \dots, a_{i+1}, f, a_{i-1}, \dots, a_1)(-) = (-1)^{\dagger} \cdot  f\big( \mu_{\mathbb{A}}(a_{i-1}, \dots, a_1, -, a_r, \dots, a_{i+1}) \big),
$$
where $\dagger=\sum_{j=1}^{r}||a_j|| +||f||$;
\item For all $f \in \mathbb{D}\Hom_{\mathbb{A}}(U,V)$ and $g \in \mathbb{D}\Hom_{\mathbb{A}}(X,Y)$,  $\mu_{\triv(\mathbb{A})}(\dots,f,\dots,g,\dots)
=0$. 
\end{enumerate}

Note that if $A$ is a finite dimensional algebra, its trivial extension when considered as an algebra coincides with its trivial extension when considered as an $A_\infty$-category with one object. If $\mathbb{A}$ denotes the associated category of indecomposable projective $A$-modules, then
    $\H^0\big(\!\Tw \triv(\bA)\big) \simeq \Kb{\triv(A)}.$

\subsection{Fukaya categories of marked surfaces}

 Partially wrapped Fukaya categories of graded marked surfaces were studies by  Auroux, Bocklandt, Haiden-Katzarkov-Kontsevich, and Lekili-Polishchuk \cite{Auroux, BocklandtPuncturedSurface, HaidenKatzarkovKontsevich, LekiliPolishchukGentle}. In order to recall their construction we first need to discuss collections of arcs and line fields on surfaces. From now on we will consider only smooth surfaces and curves.

Let $\Sigma$ be a (smooth) connected, compact oriented surface with non-empty boundary $\partial \Sigma$ together with either a finite non-empty set of marked points $\marked$ on the boundary, or with a finite non-empty set of marked points $\mathscr{P}$ in the interior of $\Sigma$ (the points of $\mathscr{P}$ will be sometimes called punctures and will play the role of the vertices of the Brauer graph embedded into $\Sigma$). We will consider \textbf{arcs} connecting marked points on $\Sigma$, an arc will be called \textbf{simple}, if it has no self intersections except for maybe at its endpoints. We will also consider \textbf{arc systems} on $\Sigma$, that is collections of mutually non-homotopic simple arcs $\mathcal{A}=\{\gamma_1,\dots,\gamma_n\}$, intersecting only transversely and at the endpoints. An arc system $\mathcal{A}$ will be called \textbf{full} if $\Sigma \setminus |\mathcal{A}|$ is a disjoint union of discs and annuli, containing exactly one boundary component $B$ of $\Sigma$ and such that there are no marked points on $B$. (In the case when there are marked points only on the boundary of $\Sigma$, all boundary arcs are required to belong to $\mathcal{A}$, see  \cite{HaidenKatzarkovKontsevich}. This condition can be relaxed to allow arc systems such that each disc from $\Sigma \setminus |\mathcal{A}|$ contains at most one boundary arc not in $\mathcal{A}$.)

An example of a full arc system is a Brauer graph $\Gamma$ on the corresponding surface $\Sigma_\Gamma$: in that case there are only punctures on $\Sigma_\Gamma$, $\mathcal{A}=\Gamma$, and $\Sigma \setminus |\mathcal{A}|$ is a collection of annuli, containing exactly one boundary component of $\Sigma_\Gamma$ with no marked points. After dragging marked points to the boundary of $\Sigma_\Gamma$ as in the construction of gentle algebras we get another arc system: $\Sigma$ has marked points only on the boundary components, $\Sigma \setminus |\mathcal{A}|$ is a disjoint union of discs and annuli.

A \textbf{line field} on $\Sigma$ is a smooth section $\eta: \Sigma \setminus \marked \rightarrow \mathbb{P}(T\Sigma)$ (or $\eta: \Sigma \setminus \mathscr{P} \rightarrow \mathbb{P}(T\Sigma)$, respectively), where $\mathbb{P}(T\Sigma)$ denotes the projectivized tangent bundle of $\Sigma$. In other words at each point of $\Sigma$, which is not a marked point, we fix a tangent line, which varies continuously on $\Sigma$. The surface $\Sigma$ together with the set of marked points and the line field will be called a \textbf{graded surface}. Using the line field $\eta$ one can equip  the arcs on $\Sigma$ with the additional data of a grading. For two graded arcs on $\Sigma$ one can define graded intersections, so each intersection comes equipped with a degree in $\mathbb{Z}$.  

One of important technical tools available on a surface equipped with a line field is the \textbf{winding number function}. It associates to any loop $\gamma$ on $\Sigma$ an integer, which coincides for loops which are regular homotopic. Roughly speaking the winding number $\omega_{\eta}(\gamma)$ counts (with sign) the number of times the derivative of $\gamma$ coincides with the line field, this can be illustrated as follows:

\begin{figure}[H]
	\begin{displaymath}
	\begin{tikzpicture}
		
		\begin{scope}[scale=2.3]	

		\def\steps{8}; 
		\def\l{0.15};

		\foreach \i in {0,1,...,\steps}
		{	
	\pgfmathsetmacro{\ll}{0.7*\l+ 0.3*\l*abs(cos((1/\steps)*\i*180))};
			\draw[gray, opacity=0.7, thick] ({\ll*cos( (1/\steps)*\i*180)},{(1/\steps)*\i+\ll*sin( (1/\steps)*\i*180))})--({\ll*cos( (1/\steps)*\i*180+180)},{(1/\steps)*\i+\ll*sin( (1/\steps)*\i*180+180))});	
				
		}

				\draw[ ->] (0,-0.2)--(0,1.2) node[left]{$\gamma$};
		\end{scope}
		
	\begin{scope}[scale=2.3, shift={(1.5,0)}]	

	\def\steps{8}; 
	\def\l{0.15};

	\foreach \i in {0,1,...,\steps}
	{	
		\pgfmathsetmacro{\ll}{0.7*\l+ 0.3*\l*abs(cos((1/\steps)*\i*180))};
		\draw[gray, opacity=0.7, thick] ({\ll*cos( -(1/\steps)*\i*180)},{(1/\steps)*\i+\ll*sin( -(1/\steps)*\i*180))})--({\ll*cos(- (1/\steps)*\i*180+180)},{(1/\steps)*\i+\ll*sin( -(1/\steps)*\i*180+180))});	
		
	}
	
	\draw[->] (0,-0.2)--(0,1.2) node[left]{$\gamma$};
	\end{scope}

	\end{tikzpicture}
	\end{displaymath}
	\caption[Contributions of a segment of $\gamma$ to its winding number with respect to a line field]{Contributions of a segment of $\gamma$ to its winding number with respect to a line field   \! \begin{tikzpicture}[scale=2, baseline=-3.3ex]\begin{scope} 
			\draw[gray, opacity=0.7, ultra thick] (0,-0.2)--(0.25,-0.2);
			\end{scope}\end{tikzpicture}. The contribution is $-1$ on the left and $+1$ on the right. Locally, we work in coordinates where $\dot \gamma(t)$ is constant.} \label{FigureWindingNumbers}
\end{figure}
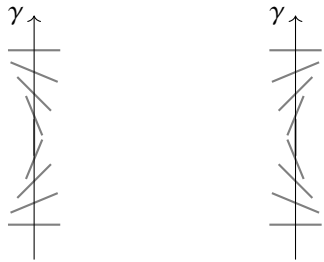

The \textbf{Fukaya category} of the graded surface $(\Sigma, \marked, \eta)$ is defined as $\Fuk(\Sigma)=\HH^0\left(\Tw \mathcal{F}\right)$ for an $A_{\infty}$-category $\mathcal{F}=\mathcal{F}_{\cA}(\Sigma)$ constructed from a full  graded arc system $\cA$. The category $\Fuk(\Sigma)$ does not depend on the choice of $\cA$ up to equivalence as shown by Haiden, Katzarkov, and Kontsevich \cite{HaidenKatzarkovKontsevich}.

 Let $D$ be a disc with $n$ marked points on the boundary. By a \textbf{marked disc} on $\Sigma$ we mean a continuous map $D \rightarrow \Sigma$ which is a smooth immersion on $D \setminus \marked$ and which sends marked points in $D$ to marked points in $\Sigma$ and boundary arcs of $D$ to arcs of $\cA$. We will assume that the edges of the disc are numbered in the clockwise order by $1,\dots,n \in \mathbb{Z}_n$. In particular, a marked disc whose corresponding sequence of arcs in $\cA$ is $\delta_1, \dots, \delta_n$ gives a cyclic sequence of \textbf{oriented intersections} $p_i \in \delta_i \oInt \delta_{i+1}$, that is locally around $p_i$ the arcs $\delta_i$ and $\delta_{i+1}$ are arranged  as in the following figure. The set of oriented intersections from $\gamma$ to $\delta$ is denoted by $\gamma \oInt \delta$. If $\alpha, \beta, \gamma$ are arcs and $p \in \alpha \oInt \beta$ and $q= \in \beta \oInt \gamma$, then   $q p \in \alpha \oInt \gamma$ will denote the corresponding intersection at the same point if it exists.

	\begin{figure}[h]
		\begin{displaymath}
		\begin{tikzpicture}[scale=0.8, hobby]
		\draw[thick, 
		](0,0) circle (50pt);
		\filldraw (-50pt, 0) circle (2pt);
		\draw (-40pt,0) node{$p_i$};
		\draw (0,1.3) node {$\partial \Sigma$};
		\draw [line width=0.5, color=black] plot  [ tension=1] coordinates {  (-50pt,0) (-3, 1) (-5, 1.5)  };
		\draw [line width=0.5, color=black] plot  [ tension=1] coordinates {  (-50pt,0) (-3, -0.5) (-5, -1.5)  };
		\draw (-3,1.3) node{$\delta_i$};
		\draw (-3,-0.8) node{$\delta_{i+1}$};
		
		\draw[dashed, ->] ({3.2*cos(165)},{3.2*sin(165)}) arc (164:188:3.2);
		\end{tikzpicture}
		\end{displaymath} 
  \caption{An oriented intersection $p_i$ from $\delta_i$ to $\delta_{i+1}$.} 
    \label{FigureIntersection}
	\end{figure}
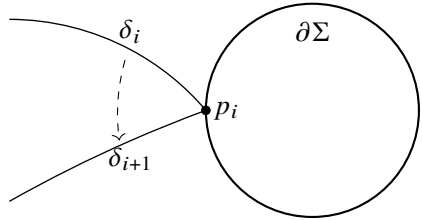

Fix a graded arc system $\cA$ on a graded surface $(\Sigma, \marked, \eta)$. The $A_{\infty}$-category $\cF=\cF_{\cA}(\Sigma)$ is an $A_{\infty}$-category with 
\begin{enumerate}[wide, labelwidth=!, labelindent=0pt]
    \item \textbf{Objects:} $\{X_\gamma\}_{\gamma\in\cA}$.
    \item \textbf{Morphisms:} For distinct $\gamma_1, \gamma_2 \in \cA$ the vector space $\Hom_{\cF}(X_{\gamma_1}, X_{\gamma_2})$ is a vector space with a basis parametrized by $\gamma_1 \oInt \gamma_2$. The identity morphism and the set of oriented self-intersections of $\gamma_1$ form a basis of $\Hom_{\cF}(X_{\gamma_1},X_{\gamma_1})$.
    \item \textbf{Grading:} The degree $|a|$ of the morphism $a$ corresponding to an intersection $p$ is given by its degree as a graded intersection of graded arcs. 
    \item \textbf{Composition:} For $a \in \Hom_{\cF}(X_{\gamma_2}, X_{\gamma_3})$ corresponding to $p \in \gamma_2 \oInt \gamma_3$ and $b \in \Hom_{\cF}(X_{\gamma_1}, X_{\gamma_2})$ corresponding to $q \in \gamma_1 \oInt \gamma_2$ put $\mu_2(a,b)=(-1)^{|b|}c$, for $c$ corresponding to $pq \in \gamma_1 \oInt \gamma_3$ if defined and $c=0$ otherwise. All other compositions of basis elements vanish. 
    \item \textbf{Higher operations:} Let $D$ be a marked disc on $\Sigma$ and denote by $p_1, \dots, p_n$ its associated sequence of oriented intersections between arcs of $\cA$. The corresponding sequence of morphisms  $a_1,\dots,a_n$  will be called  a \textbf{disc sequence}. For any disc sequence  $a_1,\dots, a_n$, $$\mu_n(a_n,\dots,a_1b)=(-1)^{|b|}b \text{ for }a_1b\neq 0 \text{ and}$$ $$\mu_n(ba_n,\dots,a_1)=b \text{ for }ba_n\neq 0.$$ 
 Higher operations vanish on all sequences of maps which are not of the form  as above.
\end{enumerate}

We can think of this $A_\infty$-category as a path category given by a quiver with relations. The arcs in $\cA$ correspond to the vertices of the quiver, the arrows are given by the ordering of arcs coming from the orientation of the surface in a similar manner to Brauer graph algebras. Forgetting the grading and the higher operations, one gets the quiver with relations of a gentle algebra. 
If there are no marked discs on $\Sigma$, all higher operations vanish, if  additionally all arrows are of degree $0$, then $\mathcal{F}_{\cA}$ is just the category of indecomposable projective modules over an ordinary gentle algebra.
Note that $\mu_1=0$ in all cases.

One can verify that $ |a_1|+\dots+|a_n|=n-2$ for any disc sequence $a_1,\dots,a_n$ which shows that all operations in $\cF_{\cA}(\Sigma)$ have the correct degree.  Note that the definition used in \cite{HaidenKatzarkovKontsevich} assumes that each boundary component of $\Sigma$ has at least one marked point. However, in the more general context considered here we assume only $\marked\neq \emptyset$ and the definition still works (see  \cite[Definition 3.24]{OpperZvonareva}).

\begin{figure}[H]
    \centering
    \begin{tikzpicture}

\node (v17) at (7.6,-2.6) {};
\node (v11) at (7.6,-3.05) {};
\draw  (v17)  ellipse (2.5 and 2.5);
\draw  (v11) ellipse (0.5 and 0.5);
\fill[gray, opacity=0.1, even odd rule]  (v17) ellipse (2.5 and 2.5) (v11) ellipse (0.5 and 0.5);

\node (v2) at (7.6,-0.1) {$\bullet$};
\node (v4) at (8.85,-4.75) {$\bullet$};
\node (v3) at (7.15,-3.3) {$\bullet$};
\draw  (7.6,-0.1) edge[bend right] (7.15,-3.3);
\draw  (8.85,-4.75) edge[bend right] (7.6,-0.1);
\draw  (7.15,-3.3) node (v9) {}  edge[bend right] (8.85,-4.75);
\node (v5) at (7.5,-4) {};
\node (v6) at (6.85,-2.2) {};
\node (v7) at (8.9,-2.3) {};

\node (v1) at (5.3,-3.5) {$\bullet$};
\draw  (5.3,-3.5) edge[bend left] (7.6,-0.1);

\draw  (5.3,-3.5)  edge (7.15,-3.3);
\draw  (5.3,-3.5) edge[bend right] (8.85,-4.75);

\node (v11) at (6.25,-3.4) {};
\node (v12)  at (5.8,-1.75) {};
\node (v13) at (6.8,-4.6) {};
\draw  (v11) edge[stealth-] (v6);
\draw  (v12) edge[-stealth] (v6);
\draw  (v11) edge[-stealth] (v12);
\draw  (v13) edge[-stealth] (v11);
\node at (6.7,-4) {\small{$b$}};
\node at (5.8,-2.8) {\small{$a_1$}};
\node at (6.45,-1.8) {\small{$a_2$}};
\node at (6.75,-2.95) {\small{$a_3$}};

\node (v200) at (5.8,-0.9) {$\bullet$};
\node (v201) at (6.8,-4.97) {$\bullet$};
\node (v202) at (10.1,-2.5) {$\bullet$};

\end{tikzpicture}

    \caption{An example of a higher operation coming from a marked disc: $\mu_n(a_3,a_2,a_1b)=(-1)^{|b|}b$.}
    \label{fig:enter-label}
\end{figure}
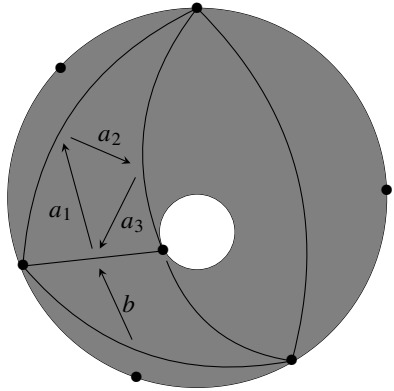

As it was mentioned before the category $\Fuk(\Sigma)=\HH^0\left(\Tw \mathcal{F}_{\cA}\right)$ does not depend on the chosen full arc system $\cA$. The proof is based on the following: For two full arc systems $\cA, \cA'$ such that $\cA'=\cA \setminus \{\gamma\}$ for some arc $\gamma \in \cA$, there is a fully faithful inclusion $\mathcal{F}_{\cA'} \subseteq \mathcal{F}_{\cA}$ of $A_{\infty}$-categories which induces a triangulated equivalence  $\HH^0\left(\Tw \mathcal{F}_{\cA'}\right) \rightarrow \HH^0\left(\Tw \mathcal{F}_{\cA}\right)$. 
Proving that the canonical inclusion $\mathcal{F}_{\cA'} \subseteq \mathcal{F}_{\cA}$ induces an equivalence uses the fact that in $\cA$ there must be a marked disc, containing the arc $\gamma$. Thus, there exists a disc sequence $a_1,\dots,a_n$ involving objects $X_1, \dots, X_n$ in $\Tw(\cF_\cA)$ with $X_1=X_{\gamma}$ and there exists a twisted complex $X_1'$ which is built from the shifts of objects $X_2,\dots,X_{n}$ and such that $X_1\simeq X_1'$ in $\H^0(\Tw\cF_\cA)$. The object $X_1'$ can be thought of as an iterated cone of objects $X_2,\dots,X_{n}$. The isomorphism $X_1\simeq X_1'$ can be constructed using the higher operation of the form $\mu_n(a_n,\dots,a_1)=1_{X_1}$. Consequently, $\cA$ and $\cA'$ generate the same triangulated subcategory of $\H^0(\Tw\cF_\cA)$ and hence the whole category $\H^0(\Tw\cF_\cA)$ (note that the objects of $\cA$ generate $\H^0(\Tw\cF_\cA)$ as a triangulated category). Having constructed the equivalences corresponding to arc inclusions, one can take any two full arc systems on $\Sigma$, complete them to triangulations of $\Sigma$, and use the fact that any two triangulations of $\Sigma$ are related by a sequence of flips (regarded as a zig-zag of two arc inclusions), proving that $\Fuk(\Sigma)$ does not depend on the full arc system $\cA$. 

The same strategy is used to show that the derived category of a Brauer graph algebra depends only on the corresponding  graded surface. 

\subsection{$A_\infty$-structures for Brauer graph algebras}

In order to carry out the plan set by the example of Fukaya categories we need to introduce some class of $A_\infty$-categories, which would contain the categories of indecomposable projectives over Brauer graph algebras. We will associate such an $A_\infty$-category to an arc system on a graded surface with punctures. In order to show that the defined structures are indeed $A_\infty$-categories we will need to put a couple of technical assumptions on the line field and on the arc system.

From now on the surface $\Sigma$ will have punctures $\mathscr{P}$ in its interior (and no marked points on the boundary). Every ribbon graph $\Gamma$ gives rise to a homotopy class of line fields on its surface $\Sigma_{\Gamma}$ as follows: For every edge $e$ of $\Gamma\subset \Sigma_{\Gamma}$ we can fix a path $\delta_e$, which crosses $e$ exactly once and transitively, does not cross any other edge of $\Gamma$ and intersects $\partial \Sigma_{\Gamma}$ transitively at each of its endpoints. The paths $\{\delta_e\}_{e\in \BE \Gamma}$ cut $\Sigma_\Gamma$ into polygons $P$ which are bounded by $\delta_e$ and segments of $\partial \Sigma$, each polygon contains exactly one vertex of $\Gamma$. These polygons can be equipped with line fields $\eta_P$ whose foliation is depicted in Figure \ref{FigureLineFieldRibbonGraph}. Since $\eta_P$ is parallel to the  paths $\delta_e$ at the boundary of $P$, the line fields $\eta_P$ can be glued to a line field $\eta_{\Gamma}: \Sigma_\Gamma \setminus \BV \Gamma \rightarrow \mathbb{P}(T\Sigma_\Gamma)$. A line field $\eta$ on a surface $\Sigma$ will be called a line field of \textbf{ribbon type} if $\Sigma$ is homotopic to $\Sigma_\Gamma$ for some ribbon graph $\Gamma$, the line field $\eta$ is homotopic to a line field $\eta_\Gamma$ constructed from $\Gamma$ as described above, and, additionally, for each $p \in \punct$ there exists a local coordinate chart around $p$ in which all leaves of $\eta$ are loops which encircle $p$. The technical advantage of considering line fields of ribbon type is that for such line fields one can define winding numbers $\omega_\eta(\gamma)$ not only for loops but also for arcs (invariantly under regular homotopy of arcs).

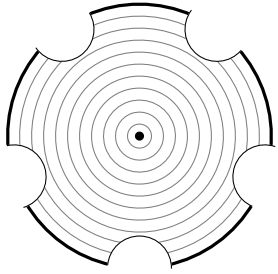
\begin{figure}[H]
\begin{tikzpicture}[scale=1.5, rotate around={18: (0,0)}]

\def\nl{10}; 
\def\r{3pt}; 
\foreach \i in {1, ..., \nl}
{
\draw[gray] (0,0) circle ({\r*\i});

}
\filldraw[black] (0,0) circle (1pt);
\def\nd{5}; 
\def\dist{\r*(\nl+1)}; 
\def\angle{360/\nd}; 
\def\length{240/\nd}; 

\foreach \i in {1,..., \nd}
{
\draw[very thick] ({\dist * cos(\angle*\i-\length*0.5)},{\dist*sin(\angle*\i-\length*0.5)}) arc ({\angle*\i-\length*0.5}:{\angle*\i+\length*0.5}:{\dist}); 
}
\def\angb{180/\nd};
\def\rtwo{8pt};
\foreach \i in {1,..., \nd}
{
\filldraw[white] ({\dist* cos(\angb * (2*\i+1)},{\dist* sin(\angb* (2*\i+1)}) circle (\rtwo);
\def\startangle{180+\angle*(\i+1)-\length*0.5-105};
\draw[black] ({\dist* cos(\angb * (2*\i+1))+\rtwo*cos(\startangle+2)},{\dist* sin(\angb* (2*\i+1))+\rtwo*sin(\startangle+2)}) arc ({\startangle+2}:{\startangle+180-2}:\rtwo);
}
\end{tikzpicture}
    \caption{The foliation of $\eta_P$ with a singularity at a vertex of $\Gamma$ in the center. Half circles in the boundary represent segments of $\partial \Sigma$. The paths $\{\delta_e\}_{e\in \BE \Gamma}$ are depicted in bold.}
    \label{FigureLineFieldRibbonGraph}
\end{figure}

A full arc system $\mathcal{A}$ on $\Sigma$ is \textbf{admissible} if every puncture $p$ can be connected to the boundary of $\Sigma$ by a path, which does not intersect the arcs of $\mathcal{A}$. 
An example of an admissible arc system on a surface with a line field of ribbon type is the collection of edges of a Brauer graph $\Gamma$ on the surface $\Sigma_\Gamma$ with the line field $\eta_\Gamma$. 

Recall that in the definition of a Brauer graph algebra as a path algebra of quiver $Q_\Gamma$ with relations, the  ideal of relations $I_\Gamma$ consisted of two types of relations: zero-relations for arrows from different cycles and relations of the form $ 
      C_\alpha^{\seq(C_\alpha)} - C_{\beta}^{\seq(C_{\beta})},$ for each edge $e$ of $\Gamma$, where $\alpha, \beta \in Q_1$ are the two arrows ending at the edge $e$. We will need to modify this definition a little bit in order to include cases with a nontrivial grading. Let $\cA$ be a full arc system on a graded surface $\Sigma$ of ribbon type and let $\seq$ be the multiplicity function on $\punct$. We can consider $\cA$ as a Brauer graph with edges corresponding to the arcs in $\cA$ (the arcs intersect only at the endpoints, so the orientation of the surface induces a cyclic ordering of arcs around the punctures of $\Sigma$). The \textbf{modified Brauer graph algebra} is the path algebra of the quiver $Q_\cA$ with relations $I'_\cA$, where the zero relations are the same as in $I_\cA$ and the second type of relations is replaced by the relations of the form $ 
      C_\alpha^{\seq(C_\alpha)} - (-1)^{\omega_\eta(\gamma)}C_{\beta}^{\seq(C_{\beta})}$ for each arc $\gamma$ of $\cA$ (here $\alpha, \beta \in Q_1$ are the two arrows ending at the edge $\gamma$). Note that if, for example, $\omega_\eta(\gamma)=0$ for all $\gamma\in\cA$ we recover the usual definition. 
Let $\delta=\alpha \pi_\mathcal{A}(\alpha)\cdots\pi_\mathcal{A}^i(\alpha)$ be a non-trivial path. We will denote by $\delta^*$ the path which complements $\delta$ to the longest path, that is $\delta\delta^*=C_\alpha^{\seq(C_\alpha)}$.
      We are finally ready to define the $A_\infty$ version of Brauer graph algebras.  

Let $\cA$ be a graded admissible arc system on a graded punctured surface $(\Sigma, \punct, \eta)$ with a line field of ribbon type. Let $\seq$ be a multiplicity function on $\punct$. The \textbf{Brauer graph category} $\mathbb{B}=\mathbb{B}(\cA, \seq)$ is the following $A_{\infty}$-category.
\begin{enumerate}[wide, labelwidth=!, labelindent=0pt]
    \item \textbf{Objects:} The set $\{X_\gamma\}_{\gamma\in \cA}$ is the set of objects of $\mathbb{B}$.
    
    \item \textbf{Morphisms:} Denote by $B'_{\mathcal{A}}$ the modified Brauer graph algebra of $(\cA, \seq)$. An arc $\gamma$ corresponds to a vertex of the quiver $Q_{\cA}$. Given graded arcs $\gamma, \delta \in \mathcal{A}$, the set of equivalence classes of paths from the vertex $\gamma$ to the vertex $\delta$ in $B'_{\mathcal{A}}$ is a $\Bbbk$-basis of $\Hom_{\mathbb{B}}(X_\gamma, X_\delta)$.
    
    \item \textbf{Grading:} Every path is a homogeneous morphism whose degree equals the sum of the degrees of its arrows. The degree $|a|$ of an arrow $a \in \Hom_{\mathbb{B}}(X_\gamma, X_\delta)$ which corresponds to an intersection $p \in \gamma \oInt \delta$ at some puncture is the degree of the intersection. The degree of the loop corresponding to any vertex of $\cA$ with valency $1$ is $0$.
    
    \item \textbf{Composition:}  For classes of paths  $b \in \Hom_{\mathbb{B}}(X,Y)$, $a \in \Hom_{\mathbb{B}}(Y,Z)$ set  $$\mu_2(a,b)\coloneqq(-1)^{|b|}a b, $$
    where $ab$ denotes the associated product in $B'_{\cA}$.
    
    \item \textbf{Higher operations:} Let $a_n,\dots,a_1$ be a disc sequence (paths $a_i$ correspond to oriented intersections $p_i\in \delta_i \oInt \delta_{i+1}$ as in the case of Fukaya categories), let $b$ be a path, then 

       \begin{equation}\label{Mu1}\mu(ba_n,\dots,a_1)=b \mbox{, for } ba_n\neq 0;\end{equation}
         \begin{equation}\label{Mu2}\mu(a_n,\dots,a_1b)=(-1)^{|b|}b \mbox{, for } a_1b\neq 0;\end{equation}
         \begin{equation}\label{HigherOps5}
         \begin{aligned} \mu(a_n,\dots, a_{r+1}, a_r(ba_r)^*,ba_r,a_{r-1},\dots, a_2)=(-1)^{\circ}a_1^* \mbox{, for } ba_r\neq 0,\\
         \text{where }\circ=|a_1|+|a_2|+\cdots+ |a_{r-1}|+|ba_r|+\omega_\eta(\delta_2)+\cdots+\omega_\eta(\delta_r).
         \end{aligned}
         \end{equation}
         \item[] Higher operations vanish for all sequences of elements which are not of the form above.
\end{enumerate}

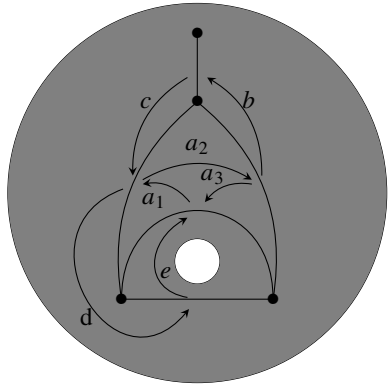
\begin{figure}[H]
    \centering
    \begin{tikzpicture}

\node (v7) at (0.5,-2.6) {};
\node (v1) at (0.5,-3.5) {};
\draw  (v7)  ellipse (2.5 and 2.5);
\draw  (v1) ellipse (0.3 and 0.3);

\fill[gray, opacity=0.1, even odd rule]  (v7)  ellipse (2.5 and 2.5) (v1) ellipse (0.3 and 0.3);

\node (v2) at (0.5,-1.4) {$\bullet$};
\node (v4) at (1.5,-4) {$\bullet$};
\node (v3) at (-0.5,-4) {$\bullet$};
\draw  (0.5,-1.4) edge[bend right] (-0.5,-4);
\draw  (1.5,-4) edge[bend right] (0.5,-1.4) ;
\draw  (-0.5,-4) edge (1.5,-4);
\node (v5) at (0.5,-4) {};
\node (v6) at (-0.35,-2.5) {};
\node (v7) at (1.35,-2.5) {};

\draw  (-0.5,-4) edge[out=90,in=90, looseness=2] (1.5,-4);

\node (v9) at (0.5,-0.5) {$\bullet$};
\draw  (0.5,-1.4) edge (0.5,-0.5);
\node (v10) at (0.5,-1) {};
\node (v8) at (0.5,-2.85) {};
\draw  (v8) edge[stealth-,bend left]node[pos=.2,above]{\small{$a_3$}} (v7);

\draw  (v6) edge[-stealth,bend left]node[pos=.5,above]{{\small{$a_2$}}} (v7);
\draw  (v6) edge[stealth-,bend left]node[pos=.2,below]{{\small{$a_1$}}} (v8);
\draw  (v7) edge[-stealth,bend right]node[pos=.5,above]{{\small{$b$}}} (v10);
\draw  (v10) edge[-stealth,bend right]node[pos=.5,above]{{\small{$c$}}} (v6);
\draw  (v6) edge[-stealth,out=200,in=-130, looseness=2.2] node[pos=.5,below]{{\small{d}}}(v5);
\draw  (v5) edge[-stealth,out=170,in=-150, looseness=1.5] node[pos=.2,above]{{\small{$e$}}} (v8);
\end{tikzpicture}
    \caption{An example of a higher operation of the form \ref{HigherOps5}: $\mu(a_3,a_2c,ba_2)=(-1)^\circ  ed$. (Here $\seq(p)=1$ for the relevant punctures $p\in \punct$.)}
    \label{fig:enter-label}
\end{figure}

One can see that the formulae \eqref{Mu1} and \eqref{Mu2} coincide with the corresponding formulae for Fukaya categories. Note that for a Brauer graph $(\Gamma,\seq)$ on the surface $\Sigma_\Gamma$ with the corresponding line field $\eta_\Gamma$ the category $\mathbb{B}(\Gamma, \seq)$ is equivalent to the category of indecomposable projectives over $B_\Gamma$, and so $\HH^0\left(\Tw \mathbb{B}(\Gamma, \seq)\right)\simeq \Kb {B_\Gamma}$. In particular, $\omega_{\eta_\Gamma}(e)=0$ for all edges of $\Gamma$ and all arrows have degree zero.

The fact that these formulae yield a structure of an $A_\infty$-category can be checked in two steps. First, when all multiplicities of the punctures are $1$,  $\mathbb{B}(\cA, \seq)$ is isomorphic to the $A_\infty$ trivial extension of the $A_\infty$-category $\mathcal{F}_{\cA}$ of the graded surface $\Sigma'$ obtained from $\Sigma$ by dragging the punctures to the boundary of $\Sigma$ (and thus turning them into marked points on the boundary) along some paths connecting the punctures to the boundary of $\Sigma$. This is possible, since $\cA$ is admissible. The arc system on $\Sigma$ can be considered as an arc system on $\Sigma'$ and thus the $A_\infty$-category $\mathcal{F}_{\cA}$ is well defined.  

The second step is the $A_\infty$ version of the construction from \cite{GreenSchrollSnashall,Asashiba}. For an admissible arc system $\cA$ on a surface of ribbon type $(\Sigma, \punct, \eta)$ with non-trivial multiplicities one can construct a branched cover $\tilde{\Sigma}\xrightarrow{\Theta} \Sigma$ ramified over $\punct$, a system of arcs $\tilde{\cA}$ on $\tilde{\Sigma}$ and a line filed $\tilde{\eta}$ equivariant under the action of the deck transformation group $\Deck(\Theta)$ of the covering. To the arc system $\tilde{\cA}$ on $(\tilde{\Sigma},\tilde{\punct}, \tilde{\eta})$ with multiplicities of punctures in  $\tilde{\punct}$ equal $1$ one can associate an $A_\infty$ Brauer graph category $\mathbb{B}(\tilde{\cA}, \tilde{\seq})$. The action of $\Deck(\Theta)$ on $\tilde{\Sigma}$ restricts to an action on $\mathbb{B}(\tilde{\cA}, \tilde{\seq})$. The orbit category of $\mathbb{B}(\tilde{\cA}, \tilde{\seq})$ under this group action is an $A_\infty$-category, which contains $\mathbb{B}(\cA, \seq)$ as a subcategory (in fact they are Morita equivalent), which proves that $\mathbb{B}(\cA, \seq)$ is itself an $A_\infty$-category. 

Since the higher operations of the form $\mu_n(a_n,\dots,a_1)=1_{X_{\gamma_1}}$ are still present for Brauer graph categories in case $a_n,\dots,a_1$ is a disc sequence coming from a marked disc bounded by the arcs $\gamma_1,\dots,\gamma_n$, we can still add arcs to the arc system $\cA$ without changing the equivalence class of the category $\HH^0\left(\Tw \mathbb{B}(\cA, \seq)\right)$ as long as the arc system stays admissible. Using this and some combinatorial arguments one can prove that for two graded admissible arc systems $\cA$ and $\cB$ on a surface $\Sigma$ and a multiplicity function $\seq$ the categories $\HH^0\left(\Tw \mathbb{B}(\cA, \seq)\right)$ and $\HH^0\left(\Tw \mathbb{B}(\cB, \seq)\right)$ are equivalent.

Finally, studying the behaviour of the Brauer graph categories under homotopy of line fields and under the action of the mapping class group $\MCG(\Sigma, \punct)$ one can prove the following theorem (see \cite[Proposition 7.6]{OpperZvonareva}):

\begin{thm}
Let $\eta_1, \eta_2$ be line fields on $(\Sigma, \punct)$ of ribbon type and let $\mathcal{A}_i$ be an $\eta_i$-graded admissible arc system. If the homotopy classes of $\eta_1$ and $\eta_2$ lie in the same $\MCG(\Sigma, \punct)$-orbit, then for all multiplicity functions $\seq_1, \seq_2$ with the same multisets of multiplicities the categories $\HH^0\left(\Tw \mathbb{B}(\mathcal{A}_1, \seq_1)\right)$ and $\HH^0\left(\Tw \mathbb{B}(\mathcal{A}_2, \seq_2)\right)$ are equivalent. 
\end{thm}

In particular, the category $\HH^0\left(\Tw \mathbb{B}(\mathcal{A}, \seq)\right)$ does not depend on the choice of the admissible arc system as in the case of Fukaya categories.

\subsection{Orbits of line fields under the mapping class group action}
Specifying the last theorem to the case of Brauer graph algebras, one gets the following (see \cite[Theorem B]{OpperZvonareva}):

\begin{thm}
Let $B$ and  $B'$ be two ordinary Brauer graph algebras with Brauer graphs $\Gamma$ and $\Gamma'$ and assume that $B$ and $B'$ are not local. The algebras $B$ and $B'$ are derived equivalent if and only if their associated graded surfaces are diffeomorphic and their multisets of multiplicities coincide.
\end{thm}

The only thing left to check in order to prove Theorem \ref{IntroTheoremCriterionDerivedEquivalence} is that Brauer graph algebras with the same combinatorial invariants have diffeomorphic graded surfaces $(\Sigma_\Gamma, \BV \Gamma,\eta_\Gamma)$. Luckily, Lekili and Polishchuk  classified  $\MCG(\Sigma, \punct)$-orbits of line fields using invariants of the winding number function of the line field (in fact their results are more general allowing for marked points on the boundary of $\Sigma$) \cite{LekiliPolishchukGentle} (see \cite{Kawazumi} for the vector field case). To recall their results we need to introduce these invariants. Let $\eta$ be a line field on $\Sigma\setminus\punct$. We will denote by  $\Omega_{\eta}^{\partial}$ the multiset of winding numbers of the boundary components of $\Sigma$, that is the multiset of values $\omega_\eta(\gamma_i)$, where $\gamma_i$ are loops homotopic to the boundary component of $\Sigma$ with their induced orientations. Similarly, we will denote by $\Omega_{\eta}^{\punct}$ the multiset of winding numbers of puncture of $\Sigma$, that is  the multiset of values $\omega_\eta(\gamma_i)$, where $\gamma_i$ are loops around the punctures of $\Sigma$. We will denote by   $\gcd(\eta)$ the greatest common divisor of all integers $\omega_\eta(\gamma)$, where $\gamma$ is a non-separating simple loop on $\Sigma$. And we will denote by $\sigma(\eta)$ the invariant which takes values $0$ and $1$, and is $0$ if and only if the line field $\eta$ comes from a nowhere vanishing vector field. In terms of the winding number function $\sigma(\eta)=0$ if and only if $\omega_\eta(\gamma)$ is even for any loop $\gamma$  on $\Sigma$. The last invariant is the so-called Arf-invariant $\Arf(\eta)$, which we do not need for the classification of the Brauer graph algebras up to derived equivalence.

\begin{thm}[{\cite[Theorem 1.2.4.]{LekiliPolishchukGentle}}]\label{TheoremClassificationOrbits}
Let $\eta, \eta'$ be line fields on $\Sigma\setminus\punct$ and let $g$ denote the genus of $\Sigma$. The homotopy classes of $\eta$ and $\eta'$ are in the same $\MCG(\Sigma, \punct)$-orbit if and only if $\Omega_{\eta}^{\partial}=\Omega_{\eta'}^{\partial}$ and $\Omega_{\eta}^{\punct}=\Omega_{\eta'}^{\punct}$ as multisets and moreover one of the following conditions holds:

\begin{enumerate}
    \item $g=0$;
    \item $g=1$ and $\gcd(\eta)=\gcd(\eta')$;
    \item $g \geq 2$ and $\sigma(\eta)=\sigma(\eta')$, additionally if $\sigma(\eta)=0=\sigma(\eta')$ and $b \equiv 2 \mod 4$ for all $b \in \Omega_{\eta}^{\partial}\cup \Omega_{\eta}^{\punct}$, then $\Arf(\eta)=\Arf(\eta')$. 
\end{enumerate}
\end{thm}

The Arf-invariant does not appear, since the winding number of all punctures are $0$, by construction of the line field of ribbon type. It turns out that $\sigma(\eta_\Gamma)=0$ if and only if $\Gamma$ is bipartite. The invariant $\gcd(\eta_\Gamma)$ captures the bipartivity of $\Gamma$ as well. Namely, it is $2$ if and only if $\Gamma$ is bipatite and $1$ otherwise. The winding numbers of boundary components of $\Sigma_\Gamma$ are exactly the perimeters off the corresponding faces of $\Gamma$ up to sign, which is indicated on the following picture:

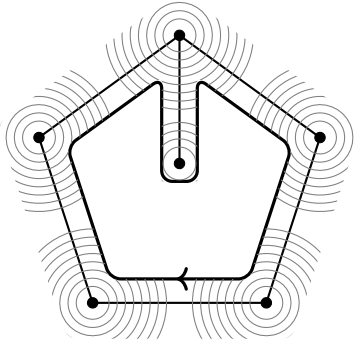
\begin{figure}[H]
  \begin{tikzpicture}
  
  \begin{scope}[rotate=18, scale=1.3]
  
  \def\r{1.5};
  \def\n{5};
  \def\angle{360 / \n};
  \def\factor{1.5};
 \draw[black, thick, opacity=100] ({1*\angle}:\r) \foreach \i in {2,..., \n} {
            -- ({\i*\angle}:\r)
        } -- cycle;
\foreach \i in {1, ..., \n} \filldraw[black] ({\i*\angle}:\r) circle (1.5pt);
    \draw [black, thick] ({\angle}:\r)--({\angle}:{0.2});
     \filldraw[black] ({\angle}:{0.2}) circle (1.5pt);

\begin{scope}        
   
\path[rounded corners, clip] ({1*\angle}:{\r*1.3}) \foreach \i in {2,..., \n} { -- ({\i*\angle}:{\r*1.3})
         } -- cycle;
         
\def\fac{0.8};
\def\dis{0.19}; 
\coordinate (D) at ({\dis-\dis*cos(72)},{-\dis*sin(72)});
\coordinate (E) at ({\dis*cos(144)-\dis*cos(72)},{\dis*sin(144)-\dis*sin(72)});

\coordinate (F) at ($ ({\angle}:{\fac*\r})+  (D) $);
\coordinate (G) at ($ ({\angle}:{0.15})+  (D) $);
\coordinate (H) at ($ ({\angle}:{0.15})+  (E) $);
\coordinate (J) at ($ ({\angle}:{\fac*\r})+  (E) $);

\draw[rounded corners, very thick, decoration={markings, mark=at position 0.2375 with {\arrow{<}}},
		postaction={decorate}] ({2*\angle}:{\r*\fac}) \foreach \i in {3,...,5} { --({\i*\angle}:{\r*\fac})} --(F)--(G)--(H)--(J)--cycle;

\begin{pgfinterruptboundingbox}
 \path[invclip] ({2*\angle}:{\r*\fac}) \foreach \i in {3,...,5} { --({\i*\angle}:{\r*\fac})}--(F)--(G)--(H)--(J)--cycle;
 \end{pgfinterruptboundingbox}

 \def\ct{9};
     \def\radius{0.75};
     \def\st{\radius / \ct};         
  \foreach \i in {1, ..., \n}
  {
     \begin{scope}[shift={({\i*\angle:\r})}]
    
     \foreach \j in {2,...,\ct}
     \draw[gray] (0,0) circle ({\st*\j});
     \end{scope}
  }      
  
  \begin{scope}

   \begin{scope}[shift={({\angle:0.2})}]
    
     \foreach \j in {2,...,5}
     \draw[gray] (0,0) circle ({\st*\j});
     \end{scope}
     
    \end{scope}

 \end{scope}

  \end{scope}
  \end{tikzpicture}
  \caption{A face of length $7$ in a ribbon graph and its line field $\eta_{\Gamma}$. Each corner of the face contributes $-1$  to the winding number of the associated boundary component, giving $-7$ as the winding number of this boundary component.}
  \label{FigureWindingNumberBoundary}
 \end{figure}

Summarising, we see that if all the combinatorial invariants of two Brauer graphs $(\Gamma_1, \seq_1)$ and $(\Gamma_2,\seq_2)$ coincide, their surfaces $\Sigma_{\Gamma_1}$ and $\Sigma_{\Gamma_2}$ are homotopic and the homotopy classes of their line fields $\eta_{\Gamma_1}$ and $\eta_{\Gamma_2}$ belong to the same orbit under the mapping class group action. Thus, the corresponding categories $\HH^0\left(\Tw \mathbb{B}(\Gamma_1, \seq_1)\right)$ and $\HH^0\left(\Tw \mathbb{B}(\Gamma_2, \seq_2)\right)$ are equivalent, and so $\Kb {B_{\Gamma_1}}$ and $\Kb {B_{\Gamma_2}}$ are equivalent as well. Finally, by \cite{RickardMoritaTheory}, the algebras are derived equivalent.


\begin{ack}
The author would like to thank Sebastian Opper for helpful comments on an earlier version of this note. She would also like to thank Mikhail Antipov and Sebastian Opper for a fruitful collaboration.
\end{ack}

\begin{funding}
This work was carried out at the Institute of Algebra and Number Theory, University of Stuttgart and at the Institute of Mathematics, Czech Academy of Sciences. It was partially supported by the  the German Research Foundation (DFG project number 512295948) and the Czech Academy of Sciences CAS  (RVO 67985840).
\end{funding}










\bibliography{bibliography}{}
	\bibliographystyle{alpha}

\end{document}